\documentclass[10pt,twoside]{article}
\usepackage{amsmath,amssymb,amscd,latexsym}
\usepackage[all]{xy}
\newcommand{\C}{{\mathbb{C}}}
\newcommand{\N}{{\mathbb{N}}}

\newcommand{\R}{{\mathbb{R}}}
\newcommand{\T}{{\mathbb{T}}}

\newcommand{\be}{\mathbf{1}}
\newcommand{\cpr}{\mathrm{cpr }\,}
\newcommand{\hr}{\mathrm{hr }\,}
\newcommand{\ev}{\mathrm{ev }}
\newcommand{\her}{\mathrm{her}}
\newcommand{\bi}{{\bar{\mbox{\it \i}}}}
\newcommand{\bj}{\bar{\mbox{\it \j}}}
\newcommand{\id}{\mathrm{id}}
\newcommand{\ord}{\mathrm{ord }\,}
\newcommand{\Prim}{\mathrm{Prim }\,}
\newcommand{\Span}{\mathrm{span}\,}

\newcommand{\tr}{\mathrm{tr}\,}

\newcommand{\Bh}{{\mathcal B}}
\newcommand{\Ch}{{\mathcal C}}
\newcommand{\Dh}{{\mathcal D}}
\newcommand{\Eh}{{\mathcal E}}

\newcommand{\Gh}{{\mathcal G}}
\newcommand{\Hh}{{\mathcal H}}
\newcommand{\Kh}{{\mathcal K}}

\newcommand{\Rh}{{\mathcal R}}

\newcommand{\Uh}{{\mathcal U}}
\newcommand{\Vh}{{\mathcal V}}
\newcommand{\Zh}{{\mathcal Z}}
\newcommand{\ot}{\overline{t}}
\newcommand{\tei}{\, | \,}

\newcommand{\verk}{\mbox{\scriptsize $\,\circ\,$}}
\newcommand{\halb}{\frac{1}{2}}
\newcommand{\te}{\textstyle}

\def\stackrellow#1#2{\mathrel{\mathop{#1}\limits_{#2}}}

\newcounter{number}[subsection]
\newcounter{altnumber}[section]
\newenvironment{nummer}{\refstepcounter{number}{\noindent\bf\arabic{section}.\arabic{subsection}.\arabic{number}}}{}
\newenvironment{altnummer}{\refstepcounter{altnumber}{\noindent\bf\arabic{section}.\arabic{altnumber}}}{}
\newcommand{\bn}{\noindent\begin{nummer} \rm}
\newcommand{\en}{\end{nummer}}

\newcommand{\altbn}{\noindent \begin{altnummer} \rm}
\newcommand{\alten}{\end{altnummer}}

\newenvironment{theorem}{\noindent {\bf Theorem:} \it}{}

\newenvironment{lemma}{\noindent {\bf Lemma:} \it}{}

\newenvironment{prop}{\noindent {\bf Proposition:} \it}{}

\newenvironment{defn}{\noindent {\bf Definition:} \it}{}

\newenvironment{cor}{\noindent {\bf Corollary:} \it}{}

\newenvironment{remark}{\noindent {\bf Remark:}}{}

\newenvironment{nremarks}{\noindent {\bf Remarks:}}{}

\newenvironment{nexamples}{\noindent {\bf Examples:}}{}
\newenvironment{proof}{\noindent {\bf Proof:}}{\mbox{}\hfill$\Box$}
\def\theequation{$\ast$}
\begin{document}

\title{{\Large \sc Covering Dimension for\\ Nuclear $C^*$-Algebras II}\vspace{3ex}}
\author{{\large Wilhelm Winter }\\ Mathematisches Institut\\ Universit\"at M\"unster\\ 
Einsteinstr. 62\\ 48149 M\"unster\\ Germany\\ e-mail: wwinter@uni-muenster.de }
\date{June 2001}
\maketitle
\setcounter{section}{-1}

\begin{abstract}
The completely positive rank is an analogue of topological covering dimension, defined for nuclear $C^*$-algebras via completely positive approximations. These may be thought of as simplicial approximations of the algebra, which leads to the concept of piecewise homogeneous maps and a notion of noncommutative simplicial complexes.\\
We introduce a technical variation of the completely positive rank and show that the two theories coincide in many important cases. Furthermore we analyze some of their properties; in particular we show that both theories behave nicely with respect to ideals and that they coincide with covering dimension of the spectrum for certain continuous trace $C^*$-algebras.
\end{abstract}

\section[{\sc Introduction}]{\large \sc Introduction}
The completely positive rank is a notion of covering dimension for nuclear $C^*$-algebras and was introduced in \cite{Wi1}. The definition is based on regarding a completely positive approximation $(F,\psi,\varphi)$ of a $C^*$-algebra $A$ as an analogue of a partition of unity. This in turn yields an analogue of an open covering of the noncommutative space $A$. The order of a covering is then modelled by a condition on the behavior of $\varphi$ on the minimal projections of $F$; this condition in some sense measures how far $\varphi$ is from preserving orthogonality. To be more precise, recall the following 

\altbn{\label{cpr}} 
\begin{defn}
Let $A,F$ be $C^*$-algebras, $F$ finite-dimensional.\\
a) We say a set $\{ e_0 , \ldots , e_n \} \subset F$ is elementary, if the $e_i$ are mutually orthogonal minimal projections.\\
b) A completely positive map $\varphi : F \to A$ is of strict order not exceeding $n$, $\ord \varphi \le n$, if the following holds:\\
For every elementary set $\{ e_0 , \ldots , e_{n+1} \} \subset F$ there are $i,j \in \{ 0 , \ldots , n+1 \}$ such that $\varphi(e_i) \perp \varphi(e_j)$.\\ 
c) $A$ has completely positive rank less than or equal to $n$, $\cpr A \le n$, if there is a net $(F_\lambda, \psi_\lambda, \varphi_\lambda)_\Lambda$ of c.p.\ approximations for $A$ such that $\varphi_\lambda \verk \psi_\lambda \to \id_A$ pointwise and $\ord \varphi_\lambda \le n \; \forall \lambda$.
\end{defn}

It then turns out that the completely positive rank has nice abstract properties; it coincides with ordinary covering dimension of the spectrum for commutative $C^*$-algebras and identifies $AF$ algebras as zero-dimensional $C^*$-algebras.
\alten

Partitions of unity may also be thought of as simplicial approximations of the space in question. The present paper is an attempt to carry this concept over to the noncommutative case in a suitable way. To this end, we have to take a closer look at  the c.p.\ approximations that compute the completely positive rank.

One might ask to what extend condition b) on the maps $\varphi$ is natural. For example, besides the minimal projections there is another class of distinguished projections in $F$, namely the (minimal) central ones. We can use these as well to define the order of $\varphi$. But then to obtain any topological invariants we have to make sure that the approximating algebra $F$ contains enough central projections to reflect some of the structure of $F$. This is done by imposing an extra condition on $\varphi$, which leads to the concept of piecewise homogeneous maps. We then obtain a variation of the completely positive rank, the homogeneous rank. This is a little easier to deal with; for example, there is an obvious (partial) result on the behavior of the homogeneous rank for tensor products. The two theories coincide in many (if not all) cases, in particular for all simple $C^*$-algebras (and for all concrete examples we have considered so far).\\
If $A$ is commutative, a partition of unity of $\hat{A}$ induces a $*$-homomorphism $\Ch(|\Sigma|) \to A$, where $|\Sigma|$ is the geometric realization of some simplicial complex $\Sigma$. These $*$-homomorphisms indeed may be thought of as simplicial approximations of $\hat{A}$. In \cite{Cu5}, Cuntz has introduced a notion of noncommutative simplicial complexes. These are universal $C^*$-algebras which (in the place of $\Ch(|\Sigma|)$) can be used to transform statements about partitions of unity (given by c.p.\ maps) into statements about $*$-homomorphisms. We apply this concept to our situation to obtain an alternative description of piecewise homogeneous maps.

Completely positive approximations with piecewise homogeneous maps are, we think, interesting by themselves since they determine something like a 'piecewise linear topology' (cf.\ \cite{BK3} for a more general approach to this idea), but the techniques developed to analyze them are also useful to compute the completely positive and the homogeneous rank in certain cases.\\ 
From the Choi--Effros lifting theorem it follows that the completely positive and the homogeneous rank behave nicely with respect to quotients. This corresponds to the (trivial) fact that $\dim K \le \dim X$ if $K$ is a closed subset of the locally compact space $X$. For an open subset $U \subset X$ we also have $\dim U \le \dim X$, but there is something to prove. In the noncommutative situation we can show that both our theories behave well with respect to ideals.\\
A $C^*$-algebra may have some obvious underlying topological space with nice properties. In this case the noncommutative covering dimension of the algebra should somehow be related to ordinary covering dimension of the underlying space. As an illustration, we examine the behavior of our theories for continuous trace algebras. It turns out that we always have $\cpr A \le \hr A \le \dim \hat{A}$ and that, under some (possibly unnecessary) extra condition, we even have equality. We will use similar (but more complicated) methods in \cite{Wi2} to analyze the completely positive rank of crossed products of manifolds with minimal diffeomorphisms.

I would like to thank J.\ Cuntz for several helpful comments, especially on Section 1.3.

%\newpage
%\input{sec1}
%\newpage
%\input{sec2}
%\newpage
%\input{sec3}
%\newpage
%\input{sec4}
%\newpage
\section[{\sc Piecewise homogeneous maps}]{\large \sc Piecewise homogeneous maps and perturbations}

\subsection[{\rm Piecewise homogeneous maps}]{\sc Piecewise homogeneous maps} 

\bn
The main ingredient in our definition of noncommutative covering dimension is the strict order of maps $\varphi : F \to A$, determined by the behavior of $\varphi$ on sets of orthogonal minimal projections.\\
We could have used sets of arbitrary orthogonal projections as well. This would have yielded a different dimension theory without some of the good abstract properties of the completely positive rank (for example, Proposition 3.7 of \cite{Wi1} would not work for this theory).\\
On the other hand, there is another class of distinguished projections in finite-dimensional $C^*$-algebras, namely the (minimal) central ones. However, one cannot expect these to reflect any of the structure of $A$ without imposing extra conditions on the c.p.\ approximations. For example, without such conditions one can always assume $F$ to be a single matrix algebra, hence containing only one central projection. Also, it does not suffice to assume $\varphi$ to be, say, completely isometric, as \cite{BK3}, Theorem 5.13 shows.\\
So we are looking for nice extra conditions on c.p.\ approximations $(F, \psi, \varphi)$ which ensure us that the central projections of $F$ reflect at least some of the structure of the approximated algebra $A$.
\en

\bn
One such condition is suggested by \cite{Wi1}, Lemma 3.14, which says that, for a c.p.\ map $\varphi : M_r \to A$ with $\ord \varphi = n$, either $n = r - 1$ or $n = 0$. If $n = r-1$, the order condition gives no information, simply because the matrix algebra is too small, but if $n = 0$, there are nice ways of describing $\varphi$ explicitly (cf.\ Proposition \ref{correspondence} and \cite{Wi1}, Proposition 4.1.1 a)).\\
In \cite{Wi1}, Proposition 3.7 we saw that $\cpr (\Ch_0(X) \otimes M_r) \le \dim X$ by constructing c.p. approximations $(F, \psi, \varphi)$, where each summand of $F$ is $M_r$, $F= \bigoplus_{i=1}^s M_r$, and the restriction of $\varphi$ to each summand has strict order zero. The order of $\varphi$ then comes from the order of a c.p.\ approximation of $\Ch_0(X)$, more precisely: $\ord \varphi = \ord \varphi \verk \iota$, where $\iota : \C^s \to F$ is the canonical unital embedding. But note that $\varphi \verk \iota = \varphi|_{\Zh(F)}$, where $\Zh(F) = \C^s$ is the center of $F$; this means that the order of $\varphi$ is determined only by its behavior on the minimal central projections of $F$.\\
This is the easiest example of a more general concept: Below we shall consider maps $\varphi$ which have strict order zero on the summands of $F$ and analyze in how far the behavior of $\varphi$ on the center of $F$ still gives information about the approximated algebra $A$.
\en

\bn{\label{d-piecewise-homogeneous}}
\begin{defn}
Let $A$, $F$ be $C^*$-algebras, $F = \bigoplus_{i=1}^s M_{r_i}$ finite-dimensional, and let $\varphi : F \to A$ be c.p.c.\\
We say $\varphi$ is piecewise homogeneous (p.h.), if $\ord \varphi_i = 0 \; \forall i=1, \ldots, s$. $\varphi$ is p.h.\ of strict order $n$, if $\varphi$ is p.h.\ and $\ord \varphi \verk \iota = n$ (where $\iota : \C^s \to F$ again is the canonical unital embedding). A c.p.\ approximation $(F, \psi, \varphi)$ is p.h.\ (of strict order $n$), if $\varphi$ is.
\end{defn}
\en

\bn{\label{r-piecewise-homogeneous}}
\begin{remark}
We will justify the first part of the preceding definition in Corollary \ref{central-cpr}, where we show that a p.h.\ map is p.h.\ of strict order $n$ if and only if it has strict order $n$ in the sense of Definition \ref{cpr}.
\end{remark}
\en

\bn{\label{minimal-perturbations}}
\begin{prop}
Let $\varphi : M_{r_1} \oplus M_{r_2} \to A$ be c.p.c.\  and let $\Uh_i \subset U(M_{r_i}), \, i=1,2$, be nonempty open sets of unitaries.\\
(i) If $e_i \in M_{r_i}, \, i=1,2$, are minimal projections, then there are nonempty open subsets $\Uh'_i \subset \Uh_i$, such that $\varphi(\be_1) \varphi(\be_2) \neq 0$ implies
\[
\varphi(u^* e_1 u) \varphi(\bar{u}^* e_2 \bar{u}) \neq 0 \; \forall u \in \Uh'_1, \; \bar{u} \in \Uh'_2 \, .
\]
(ii) If $E_i \subset M_{r_i}, \, i=1,2$, are elementary subsets, then there are nonempty open subsets $\Uh'_i \subset \Uh_i$, such that $\varphi(\be_1) \varphi(\be_2) \neq 0$ implies
\[
\varphi(u^* e u) \varphi(\bar{u}^* \bar{e} \bar{u}) \neq 0 \; \forall e \in E_1, \, \bar{e} \in E_2, \, u \in \Uh'_1, \, \bar{u} \in \Uh'_2 \, .
\]
\end{prop}

\begin{proof}
(i) By \cite{Wi1}, Lemma 1.3.8(ii), there are unitaries $u_i^{(1)}, \ldots, u_i^{(r_i)} \in \Uh_i$ for $i= 1,2$ such that $h_i := \sum^{r_i}_{j=1} u_i^{(j)*} e_i u_i^{(j)}$ is invertible in $M_{r_i}$. But $h_i$ is positive, so we have $h_i \ge \lambda_i \cdot \be_i$ for some $\lambda_i > 0$.\\
Now if $\varphi(\be_1) \varphi(\be_2) \neq 0$, we have $\varphi(h_1) \varphi(h_2) \neq 0$ and therefore 
\[
\varphi(u_1^{(j)*} e_1 u_1^{(j)}) \, \varphi(u_2^{(\bj)*} e_2 u_2^{(\bj)}) \neq 0\]
for some $j \in \{1, \ldots ,r_1\}$ and $\bj \in \{1, \ldots ,r_2\}$.\\
But then there are open neighborhoods $\Uh, \, \bar{\Uh}$ of $u_1^{(j)}$ and $u_2^{(\bj)}$, respectively, such that 
\[
\varphi(u^* e_1 u) \varphi(\bar{u}^* e_2 \bar{u}) \neq 0 \; \forall u \in \Uh, \, \bar{u} \in \bar{\Uh} \, .
\]
Now $\Uh'_1 := \Uh_1 \cap \Uh$ and $\Uh'_2 := \Uh_2 \cap \bar{\Uh}$ have the desired properties.

(ii) Suppose $E_i = \{e_i^{(1)}, \ldots, e_i^{(l_i)} \}$ for $i=1,2$ and  $\varphi(\be_1) \varphi(\be_2) \neq 0$. Apply part (i) to $e_1^{(1)}, \, e_2^{(1)}$ and $\Uh_1, \, \Uh_2$ to obtain $\Uh_1^{(1)} \subset \Uh_1$ and $\Uh_2^{(1)} \subset \Uh_2$ such that
\[
\varphi(u^* e_1^{(1)} u) \varphi(\bar{u}^* e_2^{(1)} \bar{u}) \neq 0 \; \forall u \in \Uh_1^{(1)}, \, \bar{u} \in \Uh_2^{(1)} \, .
\]
Now take $e_1^{(2)}, \, e_2^{(1)}$ and $\Uh_1^{(1)}, \, \Uh_2^{(1)}$ to obtain nonempty $\Uh_1^{(2)} \subset \Uh_1^{(1)}$ and $\Uh_2^{(2)} \subset \Uh_2^{(1)}$. Proceed inductively to obtain
\[
\emptyset \neq \Uh_{i,1} := \Uh_i^{(l_i)} \subset \Uh_i^{(l_i - 1)} \subset \ldots \subset \Uh_i^{(1)} \subset \Uh_i, \; i= 1,2 \, ,
\]
such that
\[
\varphi(u^* e_1^{(j)} u) \, \varphi(\bar{u}^* e_2^{(1) } \bar{u}) \neq 0 \; \forall u \in \Uh_{1,1}, \, \bar{u} \in \Uh_{2,1} \mbox{\, and \,} j=1, \ldots, l_1 \, .
\]
Apply the above procedure to $e_1^{(j)}$ ($j= 1, \ldots , l_1$) and $e_2^{(2)}$ with $\Uh_{i,1}$ instead of $\Uh_i$. This will yield $\Uh_{i,2} \subset \Uh_{i,1} \subset \Uh_i$ such that
\[
\varphi(u^* e_1^{(j)} u) \, \varphi(\bar{u}^* e_2^{(1)} \bar{u}) \neq 0 \; \forall u \in \Uh_{1,2}, \, \bar{u} \in \Uh_{2,2}, \, j=1, \ldots ,l_1 \, ,
\]
while we still have 
\[
\varphi(u^* e_1^{(j)} u) \, \varphi(\bar{u}^* e_2^{(2)} \bar{u}) \neq 0 \; \forall u \in \Uh_{1,2}, \, \bar{u} \in \Uh_{2,2}, \, j=1, \ldots l_1 \, .
\]
Now induction yields chains of subsets
\begin{eqnarray*}
\emptyset \neq \Uh_{i,l_2} & \subset & \Uh_{i,l_2 -1}^{(l_1)} \subset \ldots \subset \Uh_{i,2} \subset \Uh_{i,1}^{(l_1)} \subset \ldots \\
\ldots & \subset & \Uh_{i,1}^{(1)} \subset \Uh_{i,1} \subset \Uh_i^{(l_1)} \subset \ldots \subset \Uh_i^{(1)} \subset \Uh_i
\end{eqnarray*}
for $i= 1,2$. Setting $\Uh'_i := \Uh_{i,l_2}$, we see that
\[
\varphi(u^* e_1^{(j)} u) \varphi(\bar{u}^* e_2^{(j')} \bar{u}) \neq 0 \; \forall u \in \Uh'_1, \, \bar{u} \in \Uh'_2, \, j \in \{1, \ldots l_1\}, \, j' \in \{1, \ldots l_2\} \, .
\]
\end{proof}
\en
 
\bn{\label{elementary-perturbations}}
\begin{prop}
Let $A$, $F$ be $C^*$-algebras, $F = M_{r_1} \oplus \ldots \oplus M_{r_s}$ finite-dimensional and $\varphi : F \to A$ c.p.c. Furthermore, let $E_i \subset M_{r_i}$ be elementary sets and $\Uh_i \subset U(M_{r_i})$ open neighborhoods of $\be_i$ for all $i$. \\
Then there are nonempty open subsets $\Vh_i \subset \Uh_i, \; i=1, \ldots , s$, such that the following holds:\\
If $i \neq \bi \in \{1, \ldots , s\}$and $\varphi(\be_i) \varphi(\be_{\bi}) \neq 0$, we have
\[
\varphi(u^* e u) \, \varphi(\bar{u}^* \bar{e} \bar{u}) \neq 0  
\]
for all $e \in E_i$, $\bar{e} \in E_{\bi}$, $u \in \Vh_i$, $\bar{u} \in \Vh_{\bi}$.
\end{prop}

\begin{proof}
This is an iterated application of Proposition \ref{minimal-perturbations}(ii): Apply \ref{minimal-perturbations}(ii) to $E_1, \, E_2$ and $\Uh_1, \, \Uh_2$ to obtain nonempty open subsets $\Uh^{(1)}_1, \, \Uh^{(1)}_2$, such that  $\varphi(\be_1) \varphi(\be_2) \neq 0$ implies $\varphi(u^* e u) \, \varphi({u'}^* e' u') \neq 0$ for $e \in E_1$, $\bar{e} \in E_2$, $u \in \Uh^{(1)}_1$ and $\bar{u} \in \Uh^{(1)}_2$.\\
Then take $E_1, \, E_3$ and $\Uh^{(1)}_1, \, \Uh_3$ to obtain $\Uh^{(2)}_1$ and $\Uh^{(1)}_3$, respectively, afterwards $E_2, \, E_3$ and $\Uh^{(1)}_2, \, \Uh^{(1)}_3$ to obtain $\Uh^{(2)}_2$ and $\Uh^{(2)}_3$. \\
Induction yields 
\[
\emptyset \neq \Vh_i := \Uh^{(s)}_i \subset \Uh^{(s-1)}_i \subset \ldots \subset \Uh^{(0)}_i := \Uh_i
\]
for $i= 1, \ldots ,s$ such that the $\Vh_i$ have the desired property. 
\end{proof}
\en

\bn{\label{central-cpr}}
\begin{cor}
Let $A, \, F$ be as above and $\varphi : F \to A$ c.p.c., $\iota : \C^s \to F$ the canonical unital embedding. Then $\ord \varphi \verk \iota \le \ord \varphi$. If $\varphi$ is p.h., we even have equality. 
\end{cor}
\en

\bn{\label{r-piecewise-homogeneous2}}
\begin{remark}
This justifies Definition \ref{d-piecewise-homogeneous}. Furthermore, one could use Corollary \ref{central-cpr} to give a slightly simplified proof of \cite{Wi1}, Proposition 3.17.
\end{remark}
\en

\subsection[{\rm Order zero maps and stable relations}]{\sc Order zero maps and stable relations}

In this subsection we take a closer look on maps of strict order zero, the building blocks of piecewise homogeneous maps. 

\bn{\label{correspondence}}
It is well-known that, if $0 \le g \le \be$ is an element of some $C^*$-algebra $A$, then $C^*(g) \subset A$ is a quotient of $\Ch_0((0,1])$, the universal $C^*$-algebra generated by one positive element of norm $\le 1$. Equivalently, a c.p.c.\ map $\varphi : \C \to A$ induces a $*$-homomorphism $\pi : \Ch_0((0,1]) \to A$ such that $\pi (h) = \varphi(1)$ (where $h= \id_{(0,1]} \in \Ch_0((0,1])$ is the canonical generator). Revisiting \cite{Wi1}, Proposition 4.1.1, we then obtain the following matrix analogue of this correspondence (as usual, $CB := \Ch_0((0,1]) \otimes B$ denotes the cone over $B$):

\begin{prop} If $\varphi : M_r \to A$ is c.p.c.\ with $\ord \varphi =0$, then there is a unique $*$-homomorphism $\pi : CM_r \to A$ such that $\pi (h \otimes x) = \varphi(x) \; \forall x \in M_r$. Conversely, any $*$-homomorphism $\pi : CM_r \to A$ induces such a c.p.c.\ order zero map $\varphi$.
\end{prop}
\en

\bn{\label{order-zero-products}}
\begin{cor}
Let $\varphi_i : M_{r_i} \to A_i$ be c.p.c.\ maps with $\ord \varphi_i = 0$, $i=1,2$. Then the induced map $\varphi_1 \otimes \varphi_2 : M_{r_1} \otimes M_{r_2} \to A_1 \otimes A_2$ has strict order zero.
\end{cor}

\begin{proof}
The $\varphi_i$ induce $*$-homomorphisms $\pi_i : CM_{r_i} \to A_i$ such that $\pi_i (h \otimes x_i) = \varphi_i(x_i)$. Furthermore, there is a c.p.c.\ map $\varphi : M_{r_1} \otimes M_{r_2} \to CM_{r_1} \otimes CM_{r_2}$ given by $\varphi(x_1 \otimes x_2) = (h \otimes x_1) \otimes (h \otimes x_2)$. Obviously $\varphi$ has strict order zero. But then $(\pi_1 \otimes \pi_2) \verk \varphi$ also has strict order zero and
\[
(\pi_1 \otimes \pi_2) \verk \varphi (x_1 \otimes x_2) = \varphi_1(x_1) \otimes \varphi_2(x_2) = (\varphi_1 \otimes \varphi_2) (x_1 \otimes x_2) \, .
\]
\end{proof}
\en

\bn{\label{relations}}
For $i, \, j = 2, \ldots, r$ consider the relations
\def\theequation{$\Rh$}
\begin{equation}{\label{R}}
\| x_i \| \le 1 , \; x_i x_j =0, \; x_i^* x_j = \delta_{i,j} x_2^* x_2 \; ;
\end{equation}
it is well-known that 
\[
CM_r \cong C^*( x_2, \ldots, x_r \, | \, \Rh ) \, ,
\]
i.e.\  $CM_r$ is the universal $C^*$-algebra with relations ($\Rh$) (cf.\ \cite{Lo}, Table 1). As a consequence, whenever $x_2, \ldots , x_r \in A$ satisfy these relations, the assignment $e_{i-1,i} \mapsto x_i$ induces a c.p.c.\ order zero map $\varphi : M_r \to A$ and any such map is determined by a finite set of relations.\\
Furthermore, these relations are weakly stable in the sense of \cite{Lo}, Definition 4.1.1; in Lemma \ref{dominated-order-zero-maps} we will use this fact to analyze perturbations of order zero maps.
\en

\bn{\label{lifting-order-zero}}
As it turns out, the relations defining $CM_r$ are even liftable, which means that $CM_r$ is projective (\cite{Lo}, Theorem 10.2.1). We then have the following lifting result for order zero maps:

\begin{prop}
Let $J \lhd A$ be an ideal and $\varphi : M_r \to A/J$ c.p.c.\ with $\ord \varphi =0$. Then $\varphi$ has a c.p.c.\ lift $\bar{\varphi} : M_r \to A$ with $\ord \varphi =0$.
\end{prop}

\begin{proof}
By Proposition \ref{correspondence} $\varphi$ induces a $*$-homomorphism $\pi : CM_r \to A/J$, which lifts to a $*$-homomorphism $\bar{\pi} : CM_r \to A$, since $CM_r$ is projective. From $\bar{\pi}$ in turn we obtain a c.p.c.\ order zero map $\bar{\varphi} : M_r \to A$ which is easily seen to lift $\varphi$.
\end{proof}
\en

\bn{\label{dominated-order-zero-maps}}
Another reason why order zero maps are nice building blocks for c.p.\ maps is that they are stable under perturbations:

\begin{lemma}
For any $r\in \N$ and $\eta > 0$ there is $\delta > 0$ such that the following holds:\\
Let $\varphi : M_r \to A$ be c.p.c.\ with $\ord \varphi = 0$; suppose that $0 \le u,h \le 1$ in $A_+$ satisfy $\| [ u, \varphi(x) ] \| < \frac{\delta}{\|x\|} \; \forall x \in M_r$ and that $\| \varphi(\be) - h\| < \delta$. Then there is a c.p.c.\ map $\bar{\varphi} : M_r \to A$ such that
\begin{itemize}
\item[(i)]
$\| \bar{\varphi} (x) - \varphi (x) u \| < \frac{\eta}{\|x\|} \; \forall x \in M_r$
\item[(ii)] $\bar{\varphi} (M_r) \subset \overline{hAh} \cap \langle u \rangle$, where $\langle u \rangle$ is the ideal generated by $u$
\item[(iii)] $\ord \bar{\varphi} =0$.
\end{itemize}
\end{lemma}

\begin{proof}
First choose $\beta > 0$, then $\alpha >0$, then $k \in \N$ and finally $\delta > 0$. It will become clear in the course of the proof, how small (or large) these constants must be. Observe that they can be defined independently of a special choice of $\varphi$, $u$ and $h$.

Set $x_j := h^\frac{1}{k} u^\halb \varphi(e_{j-1,j}) u^\halb h^\frac{1}{k}, \; j=2, \ldots, r$, then $x_j \in B:= \overline{hAh} \cap \langle u \rangle$ and it is straightforward to check that (for any $\varphi$, $u$ and $h$ satisfying the condition of the lemma) the $x_j$ satisfy the relations (\ref{R}) of \ref{relations} within $\alpha$ if only $k$ is large and $\delta$ is small enough. But the relations are weakly stable in the sense of \cite{Lo}, Definition 4.1.1. Therefore, if $\alpha$ is small enough, there are $y_j \in B$ satisfying ($\Rh$) exactly and such that $\| x_j - y_j \| < \beta$. Now by \ref{correspondence} and \ref{relations} the $y_j$ induce a c.p.c.\ order zero map $\bar{\varphi} : M_r \to B$ and again one checks that, if $k$ is large and $\beta$ and $\delta$ are small enough, we even have $\| \bar{\varphi} (x) - \varphi (x) u \| < \frac{\eta}{\|x\|} \; \forall x \in M_r$.
\end{proof}
\en

\bn{\label{dominated-order-n-maps}}
In the situation of the preceding lemma, one cannot expect $\bar{\varphi}$ to have order zero if $\varphi$ doesn't. Nontheless we have the following (weaker) perturbation result for maps of order $r-1$:

\begin{prop}
For any $r \in \N$ and $\eta > 0$ there is $\delta > 0$ such that the following holds:\\
Let $\varphi : M_r \to A$ be c.p.c. If $0 \le u, h \le 1$ in $A$ satisfy  $\| [u, \varphi(x)] \| < \frac{\delta}{\|x\|} \; \forall x \in M_r$ and $\|\varphi(\be_{M_r}) -h\| < \delta$, then there is a c.p.c.\ map $\bar{\varphi} : M_r \to A$ such that
\begin{itemize}
\item[(i)] $\| \bar{\varphi}(x) - \varphi(x) u \| < \frac{\eta}{\|x\|} \; \forall x \in M_r$
\item[(ii)] $\bar{\varphi}(\be_{M_r}) \subset \overline{h A h} \cap \langle u \rangle$.
\end{itemize}
\end{prop}

\begin{proof}
First choose $k \in \N$ and then $\delta > 0$, then for any $\varphi$, $u$ and $h$ set $\bar{\varphi} (\, . \,) := h^\frac{1}{k} u^\halb \varphi(\, .\,) u^\halb h^\frac{1}{k}$. It is not hard to see that $k$ and $\delta$ can be chosen independently of $\varphi$, $u$ and $h$ and that $\bar{\varphi}$ has the desired properties, if $k$ is large and $\delta$ is small enough.
\end{proof}
\en

\subsection[{\rm Noncommutative simplicial complexes}]{\sc Noncommutative simplicial complexes}

In \cite{Cu5}, Cuntz has introduced a notion of noncommutative simplicial complexes. Below we 
outline how this concept is related to p.h.\ maps. 

\bn
By definition, a (finite) simplicial complex $\Sigma$ is a set of subsets of a (finite) vertex set $V_\Sigma$ satisfying a certain coherence condition (see \cite{ES} for an introduction to simplicial complexes).\\ Let $\Ch_\Sigma^{ab,\be}$ be the universal abelian unital $C^*$-algebra with positive generators $h_\sigma$, $\sigma \in V_\Sigma$, and relations
\def\theequation{$\Gh1$}
\begin{equation}{\label{G1}}
h_{\sigma_0} \ldots h_{\sigma_n} = 0 \mbox{ if } \{ \sigma_0, \ldots , \sigma_n \} \not\in \Sigma
\end{equation}
and
\def\theequation{$\Uh1$}
\begin{equation}{\label{U1}}
\sum_{\sigma \in V_\Sigma} h_\sigma = \be \, .
\end{equation}
Then $\Ch_\Sigma^{ab,\be} \cong \Ch(|\Sigma|)$, where $|\Sigma|$ is the geometric realization of $\Sigma$. Changing the unitality condition to $\sum h_\sigma \le \be$, one obtains a nonunital version, $\Ch_\Sigma^{ab}$, which can be identified with $C\Ch(|\Sigma|)$.

Cuntz then defines the noncommutative simplicial complex associated to $\Sigma$, $\Ch_\Sigma^\be$, as the universal (noncommutative) unital $C^*$-algebra with positive generators $h_\sigma$, $\sigma \in V_\Sigma$, and relations (\ref{G1}) and (\ref{U1}). Note that the nonunital version $\Ch_\Sigma$ is no longer the suspension of $\Ch_\Sigma^{\be}$.
\en

\bn
One can often use universal $C^*$-algebras to transform classes of c.p.c.\ maps into $*$-homomorphisms. We have already seen an application of this concept in Proposition \ref{correspondence}. In the context of noncommutative simplicial complexes, for any $C^*$-algebra $A$ there is a natural one-to-one correspondence between $*$-homomorphisms $\Ch_\Sigma \to A$ and c.p.\ maps $\C^\Sigma \to A$ respecting (\ref{G1}).
\en

\bn
It is a remarkable fact that the natural surjection $\Ch_\Sigma^{\be} \to \Ch_\Sigma^{ab,\be}$ is a $KK$-equivalence (\cite{Cu5}, Theorem 2.13); so from this point of view, $\Ch_\Sigma^{\be}$ is the 'right' noncommutative analogue of $\Ch(|\Sigma|)$. However, the relations (\ref{G1}) are often hard to control in concrete applications. One might therefore ask if there are relations which describe intersections of open subsets in the commutative case (so do (\ref{G1})) and at the same time are tractable in the general $C^*$-algebraic context. 

Consider the relations
\def\theequation{$\Gh2$}
\begin{equation}{\label{G2}}
h_{\sigma_0} h_{\sigma_1} = 0 \mbox{ if } \{ \sigma_0, \sigma_1 \} \in \Sigma \, ;
\end{equation}
these are easier to deal with than (\ref{G1}), since they only involve orthogonality of positive elements instead of products with more than two factors. Again one can define universal $C^*$-algebras $\Dh_\Sigma^{ab,\be}$, $\Dh_\Sigma^{\be}$, $\Dh_\Sigma^{ab}$ and $\Dh_\Sigma$; note that $\Dh_\Sigma^{\be}$ is called $\Ch_\Sigma^{flag}$ in \cite{Cu5}.

In (\ref{G2}) only the 1-simplexes of $\Sigma$ occur, so we cannot expect $\Dh_\Sigma$ or $\Dh_\Sigma^{\be}$ to be good noncommutative analogues of simplicial complexes unless the structure of $\Sigma$ is already given by its 1-skeleton. $\Sigma$ is said to be full, or a flag complex, if it satisfies the following condition:
\[
\forall \, S \subset V_\Sigma, \; S \in \Sigma \Longleftrightarrow (\{s,t\} \in \Sigma \; \forall \{s,t\} \subset S ) \, .
\]
This condition implies that indeed $\Sigma$ is determined by its 1-simplexes and it turns out that $\Dh_\Sigma^{ab,\be} \cong \Ch_\Sigma^{ab,\be} (\cong \Ch(|\Sigma|))$ iff $\Sigma$ is full. From the topological point of view, we can always restrict to flag complexes, since the barycentric subdivision of any simplicial complex automatically is a flag complex. 
\en

\bn
We should mention that sometimes it is convenient to replace (\ref{U1}) by a relation like
\def\theequation{$\Uh2$}
\begin{equation}{\label{U2}}
\| \sum_{\sigma \in V_\Sigma} h_\sigma - \be \| < \varepsilon
\end{equation}
for some $\varepsilon > 0$. This is because we mainly use noncommutative simplicial complexes to describe approximations of a $C^*$-algebra; in this context one can often assume that (\ref{U2}) holds but there is no canonical way to make this relation exact without affecting (\ref{G2}).
\en

\bn
Again we have a bijection between $*$-homomorphisms $\Dh_\Sigma \to A$ and c.p.\ maps $\C^\Sigma \to A$ that respect (\ref{G2}). Furthermore, whenever $\varphi : \C^s \to A$ is c.p.c., we can associate to $\varphi$ a flag complex $\Sigma$ such that $\varphi$ respects the relations (\ref{G2}) for $\Sigma$ as follows: Define $V_\Sigma := \{1, \ldots, s \}$ and let $\{i_0,i_1\}$ be in $\Sigma$ whenever $\varphi(e_{i_0}) \varphi(e_{i_1}) \neq 0$. Then let $\Sigma$ be the flag complex generated by this set of 1-simplexes. One checks that the strict order of $\varphi$ equals the combinatorial dimension of $\Sigma$.
\en

\bn
Next let $\varphi : \bigoplus^s_{i=1} M_{r_i} \to A$ be (c.p.c.\ and) piecewise homogeneous. Of course we can apply the above procedure to the map $\varphi \verk \iota : \C^s \to A$ to obtain a flag complex $\Sigma$ and a $*$-homomorphism $\Dh_\Sigma \to A$ that is induced by $\varphi$. But it is also interesting to ask wether there are universal $C^*$-algebras for which we obtain bijections between $*$-homomorphisms and p.h.\ maps.\\
Starting with a flag complex $\Sigma$ and $F = \bigoplus_{\sigma \in V_\Sigma} M_{r_\sigma}$ we can define another version of a noncommutative simplicial complex by putting together the relations (\ref{R}) of \ref{relations} and (\ref{G2}).\\
More precisely, let $\Eh_{\Sigma,F}$ be the universal $C^*$-algebra with generators $h_{\sigma,i}$ for $\sigma \in V_\sigma$ and $i \in \{2,\ldots,r_\sigma\}$, satisfying the following relations:
\def\theequation{$\Gh3$}
\begin{equation}{\label{G3}}
\left.
\begin{array}{l}  
 h_{\sigma,i}^* h_{\sigma,j} = \delta_{i,j} \cdot h_{\sigma,2}^* h_{\sigma,2} \mbox{ for } \sigma \in V_\Sigma \mbox{ and } i,j \in \{2,\ldots,r_\sigma\} \, ,\\
\mbox{(\ref{G2})} \mbox{ with } h_\sigma := \sum_{i=2}^{r_\sigma} h_{\sigma,i} h_{\sigma,i}^* + h_{\sigma,2}^* h_{\sigma,2} \, , \\
\sum_{\sigma \in V_\Sigma} h_\sigma \le \be 
\end{array} 
\right\}
\end{equation}
(the necessary modifications if there are $\sigma$ with $r_\sigma =1$ are obvious). Then one checks that we have a bijection between p.h.\ maps and $*$-homomorphisms and that $\Dh_\Sigma \cong C^*(h_\sigma | \, \sigma \in V_\Sigma) \subset \Eh_{\Sigma,F}$.
\en

\bn
So there are various types of noncommutative simplicial complexes; all of these transform c.p.c.\ maps into $*$-homomorphisms naturally. A very interesting application of this concept is Cuntz's approach to the Baum--Connes conjecture in \cite{Cu5}, where he uses the algebras $\Ch_\Sigma^\be$ to give a conceptual explanation for the choice of the left hand side of this conjecture. \\
In our context the algebras $\Eh_{\Sigma,F}$ seem to be the suitable analogues of simplicial complexes. However, one has to be aware of the fact that the $\Eh_{\Sigma,F}$ themselves are not accessible to our theory, since in general they are far from being nuclear (and so are the $\Ch_\Sigma$ and the $\Dh_\Sigma$). 
\en

\section[{\sc Ideals}]{\large \sc Ideals}

As an application of the methods developed in the preceding section, we determine the behavior of the completely positive rank for ideals. First we need another technical observation. 

\altbn{\label{dominated-order}}
\begin{prop}
Let $A, \, F$ be $C^*$-algebras, $F = M_{r_1} \oplus \ldots \oplus M_{r_s}$ finite-dimensional, and $\varphi, \, \bar{\varphi} : F \to A$ c.p.c.\ such that
\begin{itemize}
\item[(i)] $\bar{\varphi} (\be_{M_{r_i}}) \bar{\varphi} (\be_{M_{r_\bi}}) \neq 0 \, \Longrightarrow \, \varphi(\be_{M_{r_i}}) \varphi(\be_{M_{r_\bi}}) \neq 0 \; \forall i, \bi \in \{1, \ldots , s\}$,
\item[(ii)] $\ord \bar{\varphi}_i = 0$ if $\ord \varphi_i = 0$.
\end{itemize}
Then $\ord \bar{\varphi} \le \ord \varphi$. 
\end{prop}

\begin{proof}
Let $n:= \ord \varphi$ and suppose $\ord \bar{\varphi} > n$, i.e.\ there is an elementary set $E \subset F$ with $|E| = n+2$ and 
\[
\bar{\varphi} (e) \bar{\varphi}(\bar{e}) \neq 0 \; \forall e, \bar{e} \in E \, .
\]
Set $E_i := E \cap M_{r_i}$.\\
If $|E_i| > 1$, we have $\ord \bar{\varphi}_i > 0$, hence $\ord \varphi_i > 0$ by (ii). But then by \cite{Wi1}, Lemma 3.14, we have $\ord \varphi_i = r_i - 1$, so in particular there is an elementary set $E'_i \subset M_{r_i}$ with $|E'_i| = |E_i|$ and
\[
\varphi_i(f) \varphi_i(\bar{f}) \neq 0 \; \forall f, \bar{f} \in E'_i \, .
\]
If $|E_i| =1$, let $E'_i := E_i \subset M_{r_i}$. \\
Define $I := \{ i \in \{1, \ldots ,s\} \, | \, |E_i| > 0 \}$ and note that $\sum_{i \in I} |E_i| = n+2$. For each $i \in I$, take an open neighborhood $\Uh_i \subset U(M_{r_i})$ of $\be_{M_{r_i}}$, such that
\def\theequation{$\ast$}
\begin{eqnarray}
\varphi_i(u^* f u) \, \varphi_i(u^* \bar{f} u) \neq 0 \; \forall f, \bar{f} \in E'_i, \, u \in \Uh_i \, . \label{1}
\end{eqnarray}
We have $\bar{\varphi}(\be_{r_i}) \, \bar{\varphi}(\be_{r_{\bi}}) \neq 0$ for all $i, \, \bi \in I$.\\
Let $F_I := \bigoplus_{i \in I} M_{r_i} \subset F$ and apply Proposition \ref{elementary-perturbations} to $F_I$ and $\varphi|_{F_I}$. This yields nonempty open subsets $\Vh_i \subset \Uh_i$, $i \in I$, such that the following holds:\\
For $i \neq \bi \in I$ and $f \in E'_i$, $\bar{f} \in E'_{\bi}$, $u \in \Vh_i$ and $\bar{u} \in \Vh_{\bi}$ we have
\def\theequation{$\ast \ast$}
\begin{eqnarray}
\varphi(u^* f u) \, \varphi(\bar{u}^* \bar{f} \bar{u}) \neq 0 \, . \label{2}
\end{eqnarray}
\def\theequation{$\ast$}
For each $i \in I$ choose some $u_i \in \Vh_i$ and set $E''_i := u_i^* E'_i u_i$. Then obviously each $E''_i \subset M_{r_i}$ is elementary and so is $E'' := \bigcup_{i \in I} E''_i$.\\
Using (\ref{1}) and (\ref{2}) we obtain that
\[
\varphi(e) \, \varphi(\bar{e}) \neq 0 \; \forall e, \bar{e} \in E'' \, .
\]
But $|E''| = |E| = n+2$, so $\ord \varphi > n$, a contradiction.
\end{proof}

\begin{remark} 
Condition (i) is in particular fulfilled, if
\[
\bar{\varphi}(\be_{M_{r_i}}) \subset \overline{\varphi(\be_{M_{r_i}}) A \varphi(\be_{M_{r_i}})} \; \forall i \in \{1, \ldots , s\} \, .
\]
\end{remark}
\alten

\altbn{\label{ideals}}
\begin{theorem}
Let $A$ be a separable $C^*$-algebra and $J \lhd A$ an ideal. Then $\cpr J \le \cpr A$.
\end{theorem}

\begin{proof}
Let $\varepsilon > 0$ and $a_1, \ldots ,a_k \in J$, $\|a_j\| \le 1$, be given. Choose a c.p.\ approximation $(F,\psi,\varphi)$ of $A$ for $a_1, \ldots, a_k$ within $\frac{\varepsilon}{4}$ and such that  $\ord \varphi \le n := \cpr A$. We may assume $F = M_{r_1} \oplus \ldots \oplus M_{r_s}$.\\
For those $i \in \{1, \ldots, s\}$ for which $\ord \varphi_i =0$, take $\frac{\varepsilon}{4 \cdot s}$ and $r_i$ as $\eta$ and $r$, respectively, and apply Lemma \ref{dominated-order-zero-maps} to find $\delta_i > 0$ such that the assertion of \ref{dominated-order-zero-maps} holds. For values of $i$ for which $\ord \varphi_i > 0$, apply Proposition \ref{dominated-order-n-maps} to obtain $\delta_i > 0$. Set $\delta := \min \{\delta_i \tei i=1, \ldots, s\}$.\\
By \cite{Ped} 3.12.14, $J$ has a quasicentral approximate unit $(u_\lambda)_\Lambda$. We may thus choose $u \in J_+$ with $\|u\| \le 1$ such that
\[
\| u a_j - a_j \| < \frac{\varepsilon}{4}, \; j=1, \ldots,k
\]
and
\[
\| [u, \varphi_i(x_i)] \| < \delta \, \mbox{ for all } \, x_i \in M_{r_i} \, \mbox{ with } \, \|x_i\| \le 1, \, i=1, \ldots, s \, .
\]
(For the second assertion we used that $(u_\lambda)_\Lambda$ is quasicentral for $A$ and that $\{ \varphi_i(x_i) \tei x_i \in M_{r_i}, \, \|x_i\| \le 1, \, i=1, \ldots,s \}$ is compact in $A$.)\\
Now for $i=1, \ldots,s$ by our choice of $\delta$ either the assertion of Lemma \ref{dominated-order-zero-maps} (if $\ord \varphi_i =0$) or of Proposition \ref{dominated-order-n-maps} (else) holds, so there are c.p.c.\ maps $\bar{\varphi_i} : M_{r_i} \to A$, $i=1, \ldots,s$, with the following properties:
\begin{itemize}
\item[(i)] $\| \bar{\varphi}_i(x) - \varphi_i(x) u \| < \frac{\varepsilon}{4 \cdot s} \; \forall x \in M_{r_i}$ with $\|x\| < 1$
\item[(ii)] $\bar{\varphi}_i(M_{r_i}) \subset \overline{\varphi(\be_{M_{r_i}}) A \varphi(\be_{M_{r_i}})} \cap J$
\item[(iii)] $\ord \bar{\varphi}_i = 0$ if $\ord \varphi_i = 0$.
\end{itemize}
Denote by $\bar{\psi} : J \to F$ the restriction of $\psi$ to $J \lhd A$ and define $\bar{\varphi} : F \to J$ by setting $\bar{\varphi}|_{M_{r_i}} := \frac{1}{1+ \frac{\varepsilon}{4}} \cdot \bar{\varphi}_i$, $i=1, \ldots,s$. Then $\bar{\psi}$ and $\bar{\varphi}$ are c.p.c.\ (note that $\| \sum \bar{\varphi}_i(\be_{M_{r_i}}) - \sum \varphi_i (\be_{M_{r_i}}) u \| < \frac{\varepsilon}{4}$ and that $\| \sum \varphi_i (\be_{M_{r_i}}) u \| < 1$, so $\| \sum \bar{\varphi}_i \| < 1 + \frac{1}{\varepsilon}$). Furthermore, 
\begin{eqnarray*}
\| \sum_{i=1}^s \bar{\varphi}_i \bar{\psi}_i (a_j) - a_j \| & = & \| \sum_{i=1}^s \bar{\varphi}_i \psi_i (a_j) - a_j \| \\
& \le & \| \sum_{i=1}^s \bar{\varphi}_i \psi_i (a_j) - \sum_{i=1}^s \varphi_i \psi_i (a_j) u  \| \\
& & + \| \sum_{i=1}^s \varphi_i \psi_i (a_j) u - a_j u \| \\
& & + \| a_j u - a_j \| \\
& \le & \sum_{i=1}^s \| \bar{\varphi}_i (\psi_i(a_j)) - \varphi_i (\psi_i(a_j)) u \| \\
& & + \| \varphi \psi (a_j) u - a_j u \| \\
& & + \| a_j u - a_j \| \\
& < & s \cdot \frac{\varepsilon}{4 \cdot s} + \frac{\varepsilon}{4} + \frac{\varepsilon}{4} \\
& = & \frac{3}{4} \cdot \varepsilon \, . 
\end{eqnarray*}
Because $\| x - \frac{x}{1+ \alpha} \| < \alpha$ if $\|x\| < 1 + \alpha$, we have
\[
\| \bar{\varphi} \bar{\psi} (a_j) - \sum_{i=1}^s \bar{\varphi}_i \bar{\psi}_i (a_j) \| < \frac{\varepsilon}{4} \, ,
\]
therefore $(F, \bar{\psi}, \bar{\varphi})$ is a c.p.\ approximation for $a_1, \ldots, a_k$ within $\varepsilon$.\\
We have $\ord \bar{\varphi} \le n$ by Proposition \ref{dominated-order}.
\end{proof}
\alten

\altbn{\label{r-ideals}}
\begin{remark}
It follows from our construction that the approximations of $J$ can be chosen to be piecewise homogeneous, if this is true for the approximations of $A$.
\end{remark}
\alten

\section[{\sc Homogeneous vs. completely positive rank}]{\large \sc Homogeneous vs. completely positive rank}

\subsection[{\rm Homogeneous rank}]{\sc Homogeneous rank}

\bn{\label{d-hr}}
\begin{defn}
We say the homogeneous rank of $A$ is less than or equal to $n$, $\hr A \le n$, if $A$ has a system $(F_\lambda,\psi_\lambda,\varphi_\lambda)$ of piecewise homogeneous c.p.\ approximations of strict order not exceeding $n$.
\end{defn}
\en

\bn
We obviously have $\cpr A \le \hr A$ for any $C^*$-algebra $A$, whereas it is not clear if we always have equality. \\
In \cite{Wi1}, Proposition 3.5 we saw that $\cpr \Ch_0(X) \le \dim X$ for a locally compact space $X$. But to show this we used c.p.\ approximations through finite-dimensional commutative $C^*$-algebras, and such approximations automatically are piecewise homogeneous. Together with \cite{Wi1}, Proposition 3.18, it follows that $\cpr \Ch_0(X) = \hr \Ch_0(X) = \dim X$.\\
Of course $\hr A = 0$ iff $A$ is $AF$. The two theories also coincide for all other examples we have considered so far (cf.\ \cite{Wi1}). In particular the irrational rotation algebras, the Bunce--Deddens algebras and Blackadar's simple unital projectionless algebra all have homogeneous rank one; this is because, for the computation of the completely positive rank, we already used p.h.\ approximations.\\
More generally, one can prove that $\hr A = 1$ if $\cpr A = 1$. Below we will show that $\hr A = \cpr A$ whenever $A$ is simple.
\en

\bn
Just like the completely positive rank, the homogeneous rank has nice permanence properties. In particular, the proofs of \cite{Wi1}, Section 3, show that 
\begin{itemize}
\item $\hr (A \oplus B) \le \max \{\hr A, \, \hr B\}$
\item $\hr A \le \underline{\lim} \,  \hr A_n$ if $A = \lim_\to A_n$ and that
\item $\hr (A/J) \le \hr A$ if $J \lhd A$ is an ideal.
\end{itemize}
By Remark \ref{r-ideals} we also have $\hr J \le \hr A$.
\en

\bn{\label{tensor-product}}
It is not clear how the completely positive rank behaves with respect to tensor products. For the homogeneous rank we at least have the following partial result (note that the homogeneous rank takes finite values only for nuclear $C^*$-algebras, so we do not have to specify the tensor product we are working with):

\begin{prop}
Let $A$ and $B$ be $C^*$-algebras. Then
\[
\hr (A \otimes B) \le (\hr A + 1) \cdot (\hr B +1) - 1 \, .
\]
So if $B$ is $AF$, we have $\hr (A \otimes B) \le \hr A$.
\end{prop}

\begin{proof}
Let $(F_\lambda, \psi_\lambda, \varphi_\lambda)_\Lambda$ and $(F_\gamma, \psi_\gamma, \varphi_\gamma)_\Gamma$ be systems of c.p.\ approximations for $A$ and $B$, respectively, such that the $\varphi_\lambda$ are p.h.\ of strict order $\hr A$ and the $\varphi_\gamma$ are p.h.\ of strict order $\hr B$. Then $(F_\lambda \otimes F_\gamma, \psi_\lambda \otimes \psi_\gamma, \varphi_\lambda \otimes \varphi_\gamma)_{\Lambda \times \Gamma}$ is a system of c.p.\ approximations for $A \otimes B$. From Corollary \ref{order-zero-products} it follows that $\varphi_\lambda \otimes \varphi_\gamma$ is p.h.\ and it is straightforward to check that
\[
\ord (\varphi_\lambda \otimes \varphi_\gamma) \le (\ord \varphi_\lambda + 1) \cdot (\ord \varphi_\gamma +1) -1 \, .
\]
\end{proof}

\begin{remark}
Of course an estimate like
\[
\hr (A \otimes B) \le \hr A + \hr B
\]
would be much more satisfactory. However, this would certainly be hard to obtain, even if one of the factors is, say, commutative.
\end{remark}
\en

\subsection[{\rm Simple $C^*$-algebras}]{\sc Simple $C^*$-algebras}

In this section we show that at least for simple $C^*$-algebras the homogeneous and the completely positive rank coincide. The key step is \cite{Wi1}, Lemma 3.14, but we first need a structure result for simple $C^*$-algebras.

\bn{\label{r-hereditary-matrices}}
\begin{prop} (\cite{Cu3}, Proposition 1.8) Let A be simple, $0 \not= a,b \in A_+$.\\
Then there is $0 \not= y \in A$ with $yy^* \in \overline{aAa}$ and 
$y^*y \in \overline{bAb}$.
\end{prop}

\begin{remark}
As a consequence, it is straightforward to construct nonzero $x, z, z'$ with $z \in \overline{aA_+a}$, $z' \in \overline{bA_+b}$ and $x \in A$ such that  
\[ z=xz'x^*, \; z'=x^*zx, \; z'=x^*xz', \; z=xx^*z \, . \]
In particular $\overline{zAz} \cong \overline{z'Az'}$.
\end{remark}
\en

\bn{\label{l-hereditary-matrices}}
The following observation is well-known, although we could not find an explicit proof in the literature. Recall that a $C^*$-algebra is called elementary if it is isomorphic to  $\Kh(\Hh)$ for some Hilbertspace $\Hh$.
 
\begin{lemma} Let A be simple and nonelementary, $a \in A_+, \; n \in \N$.\\
Then there are pairwise orthogonal nonzero elements $ a_1, \ldots ,a_n \in 
\overline{aA_+a}$ and $x_1, \ldots ,x_{n-1} \in \overline{aAa}$ such that  
\[
a_{i+1} = x_i^*a_ix_i, \quad a_i = x_ia_{i+1}x_i^*, 
\quad i=1, \ldots ,n-1 \, . 
\]
\end{lemma}
\begin{proof}
Since $A$ is simple and nonelementary, every irreducible representation must have empty intersection with the compacts. In particular, $A$ does not contain minimal projections; it is then straightforward to construct nonzero pairwise orthogonal positive elements $f_1, \ldots, f_n \in \overline{aAa}$. Now from an inductive argument involving Remark \ref{r-hereditary-matrices} one obtains nonzero $x_i, \, z_i, \, z'_i \in \overline{aAa}$, $i=1, \ldots ,n-1$, such that \\
(i) $z_i \in \overline{z'_{i-1} A_+ z'_{i-1}}, \: z'_i \in \overline{f_{i+1}A_+f_{i+1}} \quad 
(\mbox{set } z'_0 := f_1)$\\
(ii) $z_i = x_i z'_i x_i^*, \: z'_i = x_i^* z_i x_i, \: z_i = x_i x_i^* z_i, 
\: z'_i = x_i^* x_i z'_i$.\\
Set
\[ a_i:= x_i \ldots x_{n-1} z'_{n-1} x_{n-1}^* \ldots x_i^*, 
\: i=1, \ldots ,n-1, \; a_n:= z'_{n-1}, \] 
then 
\[ a_i = x_i a_{i+1} x_i^*, \: i=1, \ldots ,n-1. \]
Induction shows that 
$a_{i+1} \in z'_i A z'_i$ for $i=1, \ldots ,n-2$:\\
Obviously $a_{n-1} = z_{n-1} \in z'_{n-2} A_+ z'_{n-2}$. 
Assume $a_{i+1} \in z'_i A z'_i$, then
\[ a_i = x_i a_{i+1} x_i^* = x_i z'_i w z'_i x_i^* = 
x_i x_i^* z_i x_i w x_i^* z_i x_i x_i^* = z_i x_i w x_i^* z_i \in 
z'_{i-1} A z'_{i-1}. \]
Therefore
\[ x_i^* a_i x_i = x_i^* x_i a_{i+1} x_i^* x_i = a_{i+1}. \]
$a_1, \ldots ,a_n \: , \; x_1, \ldots ,x_{n-1} \in \overline{fAf}$ by construction.
\end{proof}
\en

\bn
\begin{remark} The lemma can be interpreted as follows:\\
In any hereditary $C^*$-subalgebra of a simple nonelementary
$C^*$-algebra $A$ there is (for arbitrary $n \in \N$) a hereditary 
$C^*$-subalgebra of the form $M_n(B)$.
\end{remark}
\en

\bn{\label{functions}}
Recall the following notation from \cite{Wi1}, 1.3.1: For positive numbers $\alpha$ and $\varepsilon$ define continuous positive functions on $\R$
\[
f_{\alpha,\varepsilon} (t) :=  \left\{ 
    \begin{array}{cl}
0 & \mbox{for} \; t \le \alpha \\
t & \mbox{for} \; \alpha + \varepsilon \le t \\
\mbox{linear} &  \mbox{elsewhere}
\end{array} \right.
\]
and
\[
g_{\alpha , \varepsilon} (t) :=  \left\{ 
    \begin{array}{cl}
0 & \mbox{for} \; t \le \alpha \\
1 & \mbox{for} \; \alpha + \varepsilon \le t \\
\mbox{linear} & \mbox{elsewhere \,;}
\end{array} \right.
\]
let $g_{\alpha}$ denote the characteristic function of $[\alpha , \infty]$.
\en

\bn{\label{simple-hr}}
\begin{theorem}
Let $A$ be unital and simple. Then $\hr A = \cpr A$.
\end{theorem}

\begin{proof}
We show that there is a system $(F_\lambda,\psi_\lambda,\varphi_\lambda)$ of c.p.\ approximations with $\ord \varphi_\lambda \le n := \cpr A$ and such that the summands of the $F_\lambda$ are at least $(n+2) \times (n+2)$-matrices. Then \cite{Wi1}, Lemma 3.14 will imply that the $\varphi_\lambda$ have strict order zero on all of the summands of $F_\lambda$.

So let $0 < \varepsilon < 1$ and $a_1, \ldots , a_m \in A$ with $\|a_l\| \le 1$ be given. We may assume $A$ to be infinite dimensional, for otherwise $A=M_r$ for some $r$ and there is nothing to show. Then apply Lemma \ref{l-hereditary-matrices} to obtain pairwise orthogonal $b_0, \ldots , b_{n+1} \in A_+$ and $y_0, \ldots , y_{n+1} \in A$ with $\|b_j\| = \|y_j\| =1 \; \forall j$ and such that
\[
b_0 = y_j^* b_j y_j, \; j=0, \ldots , n+1 \, .
\]
Set $\delta := \frac{1}{2 (n+2)^2} (\frac{\varepsilon}{6})^4$. Since $A$ is simple, $b:= f_{1-\delta, \delta}(b_0)$ generates $A$ as an ideal; in particular there are $k \in \N$ and $c_l, \, d_l \in A$, $l= 1, \ldots, k$, such that
\[
\be_A = \sum_{l=1}^k c_l b d_l \, .
\]
But then by \cite{Cu3}, Proposition 1.10, there are $h_\lambda \in A$, $l=1, \ldots, k$, such that
\[
\be_A = \sum_{l=1}^k h_l^* b h_l
\]
and we have
\begin{eqnarray*}
(1 - \delta) \cdot \be_A & \le & (1-\delta) \cdot \sum h_l^* b h_l \\
& \le & \sum h_l^* b^\halb b_0 b^\halb h_l \\
& \le & \sum h_l^* b h_l \\
& = & \be_A \, . 
\end{eqnarray*}
For $l=1, \ldots ,k$ and $j=0, \ldots ,n+1$ define
\[
h_l^{(j)} := y_j b^\halb h_l \, ,
\]
then we obtain for all j:
\begin{eqnarray*}
(1-\delta) \cdot \be_A & \le & \sum_{l=1}^k h_l^{(j) *} b_j h_l^{(j)} \\
& \le & \sum_{l=1}^k h_l^{(j) *} h_l^{(j)} \\
& \le & \sum_{l=1}^k h_l^* b h_l\\
& = & \be_A \, ; 
\end{eqnarray*}
note that $\| h_l^{(j)} \| \le 1 \; \forall l,j$.\\
Set $ \eta := \frac{\varepsilon^8}{4 \cdot 6^{10} k^2 (n+2)^4}$ and choose a c.p.\ approximation $(F, \psi, \varphi)$ within $\eta$ for
\[
G := \{\be_A, a_i, h_l^{(j)},  h_l^{(j)*} h_l^{(j)}, b_j \, | \,i= 1, \ldots, m, \,  l=1, \ldots, k \mbox{ and } j=0, \ldots, n+1 \}
\]
such that $\ord \varphi \le n$. Again we assume $F = M_{r_1} \oplus \ldots \oplus M_{r_s}$. Now for all $l$, $j$
\begin{eqnarray*}
 h_l^{(j)*}  h_l^{(j)} - 2 \eta & \le & \varphi \psi ( h_l^{(j)*}) \varphi \psi ( h_l^{(j)}) \\
& \le & \varphi (\psi( h_l^{(j)*}) \psi( h_l^{(j)})) \\
& \le & \varphi \psi ( h_l^{(j)*}  h_l^{(j)}) \\
& \le &  h_l^{(j)*}  h_l^{(j)} + \eta \, ,
\end{eqnarray*}
so
\[
\| \varphi(\psi( h_l^{(j)*}) \psi( h_l^{(j)})) - \varphi \psi( h_l^{(j)*}) \varphi \psi( h_l^{(j)})\| < 3 \eta \, ,
\]
hence for any $y \in F$ with $\|y\| \le 1$ we have
\[
\| \varphi( y \psi( h_l^{(j)})) - \varphi(y) \varphi \psi( h_l^{(j)})\| < (3 \eta)^\halb
\]
by \cite{Wi1}, Proposition 1.1.6. For any $j$ we have
\begin{eqnarray*}
\be_A & \ge & \varphi \psi (\be_A) \\
& \ge & \sum_{l=1}^k \varphi \psi ( h_l^{(j)*}  h_l^{(j)}) \\
& \ge & \sum \varphi (\psi( h_l^{(j)*}) \psi( h_l^{(j)})) \\
& \ge & \sum \varphi ( \psi( h_l^{(j)*}) \psi(b_j) \psi( h_l^{(j)})) \\
& \ge & \sum (\varphi \psi( h_l^{(j)*}) \varphi \psi(b_j) \varphi \psi( h_l^{(j)}) - 2 \cdot (3 \eta)^\halb \cdot \be_A) \\
& \ge & \sum ( h_l^{(j)*} b_j  h_l^{(j)} - (2 \dot (3 \eta)^\halb + 3 \eta) \cdot \be_A) \\
& \ge & (1 - \delta - k (2 (3 \eta)^\halb + 3 \eta)) \cdot \be_A \, .
\end{eqnarray*}
As a consequence,
\begin{eqnarray*}
0 & \le & \| \varphi (\be_F - \sum^k_{l=1} \psi(h_l^{(j)*}) \psi(b_j) \psi(h_l^{(j)}))\| \\
& \le & \be_A - \sum^k_{l=1} \varphi(\psi(h_l^{(j)*}) \psi(b_j) \psi(h_l^{(j)})) \| \\
& \le & \delta + k (2 \cdot (3 \eta)^\halb + 3 \eta) \\
& =: & \alpha \, .
\end{eqnarray*}
Set $q^{(j)} := g_{1-\alpha^\halb} (\sum^k_{l=1} \psi(h_l^{(j)*}) \psi(b_j) \psi(h_l^{(j)}))$, then
\[
\be_F - q^{(j)} \le \frac{1}{\alpha^\halb} (\be_F - \sum^k_{l=1} \psi(h_l^{(j)*}) \psi(b_j) \psi(h_l^{(j)}))
\]
and
\[
\varphi(\be_F - q^{(j)}) \le \frac{\alpha}{\alpha^\halb} \cdot \be_A = \alpha^\halb \cdot \be_A \,.
\]
Denote by $q_i^{(j)} := q^{(j)} \be_{M_{r_i}}$ the part of $q^{(j)}$ in $M_{r_i}$. Now if $ \| \psi_i(b_j)\| \le 1 - 2 \alpha^\halb$ for some $i \in \{1, \ldots, s\}$, then
\begin{eqnarray*}
(1 - \alpha^\halb) \cdot q_i^{(j)} & \le & q_i^{(j)}  (\sum^k_{l=1} \psi(h_l^{(j)*}) \psi(b_j) \psi(h_l^{(j)})) q_i^{(j)} \\
& \le & (1- 2 \alpha^\halb) q_i^{(j)} (\sum^k_{l=1} \psi(h_l^{(j)*}) \psi(h_l^{(j)})) q_i^{(j)} \\
& \le & (1- 2 \alpha^\halb) q_i^{(j)} (\sum^k_{l=1} \psi(h_l^{(j)*}h_l^{(j)})) q_i^{(j)} \\
& \le & (1- 2 \alpha^\halb) q_i^{(j)} \, ,
\end{eqnarray*}
which implies $q_i^{(j)} = 0$.
Let $p$ be the unit of 
\[
\bigoplus_{ \{ i \, | \, q_i^{(j)} \neq 0, \, j=0, \ldots , n+1 \} } M_{r_i} \, .
\]
Then for each $i$ we have $p \be_{M_{r_i}} \neq 0$ iff $q_i^{(j)} \neq 0 \; \forall j=0,\ldots,n+1$. Note that $\be_F - p \le  \be_F - \sum^{n+1}_{j=0} q^{(j)}$ by construction. Define $F' := pFp$, $\psi' := \psi_p$ and $\varphi' := \varphi|_{F'}$. For any $a \in A$ with $\|a\| \le 1$ we then have
\begin{eqnarray*}
\| \varphi' \psi'(a) - \varphi \psi(a) \| & = & \| \varphi(p \psi(a) p) - \varphi \psi(a) \| \\
& \le & 2 \cdot \| \varphi((\be_F - p) \psi(a) p) \| + \| \varphi((\be_F - p) \psi(a) (\be_F - p))\| \\
& < & 3 \cdot \| \varphi(\be_F - p) \|^\halb \\
& \le & 3 \cdot \|  \varphi(\be_F - \sum^{n+1}_{j=0} q^{(j)}) \|^\halb \\
& \le & 3 \cdot ((n+2) \alpha^\halb)^\halb \\
& < & \frac{\varepsilon}{2} \, ,
\end{eqnarray*}
so $(F', \psi', \varphi')$ is a c.p.\ approximation for $G$ within $\frac{\varepsilon}{2} + \eta$, which is less than $\varepsilon$. Obviously we still have $\ord \varphi' \le n$.

$F'$ consists of those summands $M_{r_i}$ of $F$ for which $q_i^{(j)} \neq 0$ for $j=0,\ldots,n+1$. This implies that (for those $i$) $\|\psi_i(b_j)\| > 1 - 2 \alpha^\halb \; \forall j$ as we have seen above. On the other hand, $\| \psi_i( \sum_{j=0}^{n+1} b_j)\| \le 1$ and we obtain
\[
1 \ge \tr (\psi_i(\sum_{j=0}^{n+1} b_j)) > (n+2) \frac{1- 2 \alpha^\halb}{r_i} \, ,
\]
where $\tr$ is the normalized trace. But then
\[
r_i > (n+2) (1 - 2 \alpha^\halb) > n+1 \, ,
\]
since $\alpha^\halb < \frac{1}{2 (n+2)}$. Applying \cite{Wi1}, Lemma 3.14, we now see that $\ord(\varphi|_{M_{r_i}}) = 0$ for all indices $i$ for which $M_{r_i}$ is a summand of $F'$, and this means that $\varphi'$ indeed is piecewise homogeneous, which was to be shown.
\end{proof}
\en

\section[{\sc Continuous trace $C^*$-algebras}]{\large \sc Continuous trace $C^*$-algebras}

In this section we examine the completely positive rank of a separable continuous trace $C^*$-algebra $A$. It turns out that $\cpr A \le \dim \hat{A}$, $\hat{A}$ being the spectrum of $A$. Under some (possibly unnecessary) extra condition, basically saying that all irreducible representations of $A$ must be of (at least locally) bounded finite dimension, we show that $\cpr A = \dim \hat{A}$.

We start by recalling some facts about the spectrum and the primitive ideal space and about $C^*$-algebras with Hausdorff spectrum and continuous trace $C^*$-algebras. Most of this material is taken from \cite{RW}; cf. also \cite{Dx} or \cite{Fe}.

\subsection[{\rm $C^*$-algebras with Hausdorff spectrum}]{\sc $C^*$-algebras with Hausdorff spectrum}

\bn
Recall that, for a $C^*$-algebra $A$, $\Prim A$ denotes the primitive ideal space, the space of kernels of irreducible representations endowed with the Jacobson topology. The spectrum $\hat{A}$ of $A$ is the space of unitary equivalence classes of irreducible representations; it inherits its topology from $\Prim A$ via the natural surjection.
\en

\bn
For a closed set $M \subset \Prim A$ we define $J_M := \bigcap \{ t \in \Prim A \tei t \in M \}$, which is a closed ideal in $A$. It is well-known that, if $\hat{A}$ is Hausdorff, then the natural map $\hat{A} \to \Prim A$ is a homeomorphism. It is then often more convenient to consider $A / J_t$ for $t \in \Prim A$ than $\pi_t (A)$, $\pi_t$ being the representation in $\hat{A}$ corresponding to $t$. We write $a (t)$ for the image of $a$ in $A / J_t$.
\en

\bn
For any $C^*$-algebra $A$ the Dauns--Hofmann Theorem (\cite{RW}, Theorem A.34) identifies $\Ch_b (\hat{A})$ with $ZM (A)$, the center of the multiplier algebra of $A$. By restriction, this makes $A$ a $\Ch_0 (\hat{A})$-module. If $\hat{A}$ is Hausdorff, the module structure is given by the formula $(f \cdot a) (t) = f (t) \cdot a (t)$ for $a \in A , f \in \Ch_0 (\hat{A}) , t \in \hat{A}$; it is easy to see that the $\Ch_0 (\hat{A})$-action is nondegenerate, i.e.\ 
\[
\overline{\Ch_0 (\hat{A}) \cdot A} = \overline{\Span \{ f \cdot a \tei f \in \Ch_0 (\hat{A}) , a \in A \}} = A \, .
\]
\en

\bn
The next result basically is a corollary of the Dauns--Hofmann Theorem (part (i) is \cite{RW}, Corollary 5.11, part (ii) is straightforward):

\begin{lemma}
  \label{nondegenerate}
Let $A$ be a $C^*$-algebra with Hausdorff spectrum, $M \subset \hat{A}$ a closed subset, $U \subset M$ open in $\hat{A}$. Then\\
(i) $J_M = \overline{\Ch_0 (\hat{A} \backslash M) \cdot A}$;\\
(ii) the quotient map $\pi_M : A \to A / J_M$ is an isomorphism on $\overline{\Ch_0 (U) \cdot A} \subset A$, which we may therefore consider as hereditary subalgebra (or even as an ideal) of $A / J_M$. (Note that we consider $\Ch_0 (\hat{A} \backslash M)$ and $\Ch_0 (U)$ as subalgebras of $\Ch_0 (\hat{A})$.)
\end{lemma}
\en

\bn
We now turn to the definition of continuous trace $C^*$-algebras. Suppose $A$ is a $C^*$-algebra with Hausdorff spectrum. We have seen that, for each $t \in \hat{A}$, $A / J_t$ has a unique irreducible representation $\pi_t$ (up to unitary equivalence); so if $p (t) \in A / J_t$ is a projection, we may define the rank of $p (t)$ as the rank of $\pi_t (p (t)) \in \Bh (\Hh_{\pi_t})$. This is well-defined since the dimension of a subspace is preserved by unitaries. In particular, we may speak about rank-one projections in $A / J_t$.
\en

\bn{\label{continuous-trace}}
\begin{defn}
A $C^*$-algebra $A$ is said to have continuous trace if it has Hausdorff spectrum and satisfies Fell's condition, i.e.\ for each $t \in \hat{A}$ there are a neighborhood $U$ of $t$ and $a \in A$ such that $a (s)$ is a rank-one projection for each $s \in U$.
\end{defn}
\en

\bn{\label{p-continuous-trace}}
There is another characterization of continuous trace algebras which shall be useful in the sequel:

\begin{prop} {\rm (\cite{RW}, Proposition 5.15)} 
Let $A$ be a separable $C^*$-algebra with Hausdorff spectrum. Then $A$ is a continuous trace $C^*$-algebra if and only if it is locally stably isomorphic to $\Ch_0 (\hat{A})$, i.e.\ each $t \in \hat{A}$ has a compact neighborhood $M$ such that $(A / J_M ) \otimes \Kh \cong \Ch (M) \otimes \Kh$ (as $C^*$-algebra and as $\Ch (M)$-modules).
\end{prop}
\en

\subsection[{\rm {\sc cpr} of continuous trace $C^*$-algebras}]{\sc Completely positive rank of continuous trace $C^*$-algebras}

\bn{\label{continuous-trace-le}}
\begin{theorem}
Let $A$ be a separable continuous trace $C^*$-algebra. Then $\hr A \le \dim \hat{A}$. 
\end{theorem}

\begin{proof}
  The argument is a generalization of Proposition \cite{Wi1}, Proposition 3.7. We use the fact that any $t \in \hat{A}$ has a compact neighborhood $K_t$ such that, locally, $A$ looks like a hereditary subalgebra of $\Ch (K_t) \otimes \Kh$, and such an algebra may be approximated by a subalgebra of $\Ch (K_t) \otimes M_r$; if $K_t$ is taken to be small enough it even suffices to consider constant functions from $K_t$ to $M_r$. This construction is used to get an open covering $(V_i)$ of strict order less than or equal to $\dim \hat{A}$, such that $(\Ch_0 (V_i) \otimes M_{r_i})$ approximates the local structure of $A$ sufficiently well. It is then easy to define $F$ and $\psi : A \to F$. To define $\varphi : F \to A$ one has to use a partion of unity to glue the $M_{r_i}$'s together. This $\varphi$ will be of strict order not exceeding $\dim \hat{A}$. Again, there is some extra work to do because we do not assume $\hat{A}$ to be compact.

All this we will now make precise. 

Consider a finite set $G \subset A_+$ with $ \| a \| \le 1$ for each $a \in G$ and $\varepsilon > 0$. We are looking for a c.p.\ approximation $(F , \psi , \varphi)$ for $G$ within $\varepsilon$ such that $\ord \varphi \le n := \dim \hat{A}$.

By Proposition \ref{p-continuous-trace} every $t \in \hat{A}$ has a compact neighborhood $K_t$ such that $(A / J_{K_t} )\otimes \Kh \cong \Ch (K_t) \otimes \Kh$ as $\Ch (K_t)$-modules. There is $0 \le h_t \in A / J_{K_t} , \; \| h_t \| \le 1$, such that 
\[
\| h_t \pi_{K_t} (a) - \pi_{K_t} (a) \| < \frac{\varepsilon}{9} \; \forall\, a \in G.
\]

Consider $h_t$ and $\pi_{K_t} (a) , \, a \in G$, as functions: $K_t \to \Kh_+$ by embedding $A / J_{K_t}$ in $(A / J_{K_t} )\otimes \Kh$ as upper left corner.

It is then routine to show that there are an open neighborhood $U_t \subset K_t$ of $t$, some $\lambda_t > 0$ and a projection valued function $q_t : U_t \to \Kh$ with the following properties:
\begin{itemize}
\item[(i)] $  q_t (s) \le \lambda_t \cdot h_t (s) \; \forall s \in U_t$,
\item[(ii)] $\| q_t (s_1) - q_t (s_2) \| < \frac{\varepsilon}{36}$ and
\item[ ] $\| \pi_{K_t} (a) (s_1) - \pi_{K_t} (a) (s_2) \| < \frac{\varepsilon}{9} \; \forall \, s_1 , s_2 \in U_t, \; a \in G$,
\item[(iii)] $\| q_t (s) h_t (s) - h_t (s) \| < \frac{\varepsilon}{9} \; \forall s \in U_t$. 
\end{itemize}
It follows that
\begin{eqnarray}
 \lefteqn{ \| q_t (s) \; \pi_{K_t} (a) (s) \; q_t (s) - \pi_{K_t} (a) (s)\|}\nonumber \\
 & \le & \| q_t (s) \; \pi_{K_t} (a) (s) \; q_t (s) - q_t (s) \; \pi_{K_t} (a) (s) \; h_t (s) \; q_t (s) \| \nonumber \\
& & + \| q_t (s) \; \pi_{K_t} (a) (s) \; h_t (s) \; q_t (s)- q_t (s) \; \pi_{K_t} (a) (s) \; h_t (s) \| \nonumber \\
& & + \| q_t (s) \; \pi_{K_t} (a) (s) \; h_t (s) - q_t (s) \; \pi_{K_t} (a) (s) \| \nonumber \\
& & + \| q_t (s) \; \pi_{K_t} (a) (s) - q_t (s) \; h_t (s) \; \pi_{K_t} (a) (s) \| \nonumber \\
& & + \| q_t (s) \; h_t (s) \; \pi_{K_t} (a) (s) - h_t (s) \; \pi_{K_t} (a) (s) \| \nonumber \\
& & + \| h_t (s) \; \pi_{K_t} (a) (s) - \pi_{K_t} (a) (s) \| \nonumber \\
& < & \frac{6}{9} \cdot \varepsilon \; \forall a \in G , s \in U_t \; . \label{4}
\end{eqnarray}
The term $q_t (s) \, \pi_{K_t} (a) (s) \, q_t (s)$ suggests that $s$ must be regarded as an element of $U_t \subset K_t$, for otherwise it is just not defined. So in terms like this we often slightly misuse our notation and write $a (s)$ for $\pi_{K_t} (a) (s)$.\\
Anyhow, we will not always explicitly indicate if we work in $\Ch (K_t) \otimes \Kh$ or $A / J_{K_t} \otimes \Kh$ (or even in $A / J_{K_t} \subset (A / J_{K_t} )\otimes \Kh$), but this should be clear from the context.

Define $K := \{ t \in \hat{A} \tei \| a (t) \| \ge \frac{\varepsilon}{9}$ for some $a \in G \}$, then $K \subset \hat{A}$ is compact (cf. \cite{RW}, Lemma A.30).

There are $t'_1 , \ldots , t'_k \in \hat{A}$ such that $K \subset \bigcup^k_1 U_{t'_i}$. Set $U_0 := \hat{A} \backslash K$, then $(U_0 , U_{t'_1} , \ldots , U_{t'_k})$ is a finite open covering of $\hat{A}$, which, by \cite{Wi1}, Proposition 2.8, has a refinement $V_1 , \ldots , V_l$ of strict order not exceeding $n$. 

We may assume that $V_i \cap K \neq \emptyset$ for $i \le m$ and $V_i \cap K = \emptyset$ for $i > m$ for some $m \in \{ 1 , \ldots, l \}$. We may further assume that for each $V_i$ there is some $t_i \in V_i$ with $t_i \notin \bigcup_{j \neq i} V_j$.

Choose a partition of unity $(g_i)_{1 , \ldots , l}$ subordinate to $(V_i)_{1 , \ldots , l}$, then 
\[
\sum^m_{i=1} g_i (s) = 1 \; \forall s \in K
\]
and $g_i (t_i) = 1 \; \forall i$. Also, the $g_i$ may be viewed as elements of $\Ch_0 (V_i)$.

For each $i \le m$ there is $j_i \in \{ 1 , \ldots , k \}$ such that $V_i \subset U_{t'_{j_i}}$; let $q_i$ be the restriction of $q_{t'_{j_i}}$ to $V_i$, then $q_i \in \Ch_b (V_i, \Kh)$ is a continuous projection valued function. From now on, we write $\overline{t}_i$ for $t'_{j_i}$. Note that we have 
\[
t_i \in V_i \subset U_{\ot_i} \subset K_{\ot_i}.
\]
By (ii) and \cite{Wi1}, Proposition 1.3.7, there are partial isometries $w_i \in \Ch_b (V_i , \Kh)$ with $\| w_i (s) - q_i (s) \| < \frac{\varepsilon}{9}$, $w^*_i (s) w_i (s) = q_i (s)$, and $w_i (s) w^*_i (s) = q_i (t_i) \; \forall s \in V_i$. Note that, by construction of the $w_i , \; w_i (t_i) = q_i (t_i)$.

Now we are ready to define $(F , \psi , \varphi)$:\\
Set $F := \bigoplus^m_{i=1} M_{r_i}$, where $r_i$ is the rank of $q_i (t_i)$ in $\Kh$. Then define $\psi : A \to F$ by 
\[
\psi (a) := {\te \bigoplus\limits^m_1} \, q_i (t_i) a (t_i) q_i (t_i);
\]
here we view $a (t_i)$ (for $t_i \in V_i \subset K_{\ot_i}$) as an element of $\Kh$ via the $\Ch_0 (\hat{A})$-maps 
\[
A \xrightarrow{\pi_{K_{\ot_i}}} A / J_{K_{\ot_i}} \hookrightarrow (A/J_{K_{\ot_i}} )\otimes \Kh \cong \Ch (K_{\ot_i}) \otimes \Kh \xrightarrow{\ev_{t_i}} \Kh \; \;  \; .
\]
Also note that $M_{r_i} \cong q_i (t_i) \Kh q_i (t_i)$. Finally, define $\varphi : F \to A$ by 
\[
\varphi ( {\te \bigoplus\limits^m_1 \, x_i}) := \sum^m_1 g_i \cdot w^*_i x_i w_i.
\]

Here, $x_i \in M_{r_i} \cong q_i (t_i) \Kh q_i (t_i) \subset \Kh$, thus $w^*_i x_i w_i \in \Ch_b (V_i , \Kh)$. But then 
\[
g_i \cdot w^*_i x_i w_i \in \Ch_0 (V_i , \Kh) \cong \Ch_0 (V_i) \otimes \Kh \subset \Ch (K_{\ot_i}) \otimes \Kh \cong (A / J_{K_{\ot_i}} )\otimes \Kh.
\]
Furthermore, 
\[
g_i \cdot w^*_i x_i w_i \le \| x_i \| \cdot g_i \cdot q_i \stackrel{(i)}{\le} \| x_i \| \cdot \lambda_{\ot_i} \cdot g_i \cdot h_{\ot_i}
\]
for $x_i \in (M_{r_i})_+$. But 
\[
g_i \cdot h_{\ot_i} \in \overline{\Ch_0 (V_i) \cdot A} \stackrellow{\subset}{\her} A / J_{K_{\ot_i}} \stackrellow{\subset}{\her} A / J_{K_{\ot_i}} \otimes \Kh \, ,
\]
where the first inclusion comes from Lemma \ref{nondegenerate}. We therefore obtain 
\[
g_i \cdot w^*_i x_i w_i \in \overline{\Ch_0 (V_i) \cdot A} \subset A \, .
\]
Thus $\psi$ and $\varphi$ are well-defined and (obviously) completely positive contractions.

For $a \in G$ and $s \in \hat{A}$ we obtain
\begin{eqnarray*}
\lefteqn{  \| \varphi \verk \psi (a) (s) - a (s)\| }\nonumber\\
& = & \| \sum^m_1 g_i (s) \cdot w^*_i (s) a (s) w_i (s) - a (s) \| \\
& \le & \| \sum^m_1 g_i (s) \cdot w^*_i (s) a(s) w_i (s) - \sum^m_1 g_i (s) \cdot q_i (s) a (s) q_i (s) \| \\
& & + \| \sum^m_1 g_i (s) \cdot q_i (s) a (s) q_i (s) - \sum^m_i g_i (s) \cdot a (s) \| \\
& & + \| \sum^m_1 g_i (s) \cdot a (s) - \sum^l_1 g_i (s) \cdot a (s) \| \\
& \le & \sum^m_1 g_i (s) \| w^*_i (s) a (s) w_i (s) - q_i (s) a (s) q_i (s) \| \\
& & + \sum^m_1 g_i (s) \| q_i (s) a (s) q_i (s) - a (s) \| \\
& & + \sum^l_{m+1} g_i (s) \| a (s) \| \\
& \le & \sum^m_1 g_i (s) \cdot \frac{2}{9} \varepsilon + \sum^m_1 g_i (s) \cdot \frac{6}{9} \varepsilon + \sum^l_{m+1} g_i (s) \cdot \frac{\varepsilon}{9} \\
& \le & \frac{9}{9} \cdot \varepsilon  = \varepsilon \; .
\end{eqnarray*}
If $s \notin V_i , \; q_i (s)$ and $w_i (s)$ are not defined, but then $g_i (s) = 0$, so all summands are well-defined. The norm estimates work in the same way:\\
We only have to check 
\[
\| w^*_i (s) a (s) w_i (s) - q_i (s) a (s) q_i (s) \| < \frac{2}{9} \cdot \varepsilon
\]
and 
\[
\| q_i (s) a (s) q_i (s) - a (s) \| < \frac{6}{9} \cdot \varepsilon
\]
for $s \in V_i$. The first is true because $\| w_i (s) - q_i (s) \| < \frac{\varepsilon}{9}$ for $s \in V_i$ (and because the $a , w_i$ and $q_i$ are normed), the latter we have already checked in (\ref{4}). 

If $i > m$, then $g_i (s) \neq 0$ only if $s \notin K$, but then $\| a(s) \| < \varepsilon$, and this yields the last estimate.

So we have seen that $(F , \psi , \varphi)$ is a c.p.\ approximation of $G$ within $\varepsilon$. It only remains to show that $\ord \varphi_i = 0 \; \forall i$ and that $\ord \varphi \le n$: 

If $e_{j_1} \perp e_{j_2}$ live in the same $M_{r_i}$, then 
\[
\varphi (e_{j_1}) \varphi (e_{j_2}) = g^2_i \cdot w^*_i e_{j_1} w_i w^*_i e_{j_2} w_i = g^2_i \cdot w^*_i e_{j_1} e_{j_2} w_i = 0
\]
(recall that we embedded $M_{r_i}$ in $\Kh$ as $q_i (t_i) \Kh q_i (t_i)$ and that $w_i w^*_i = q_i (t_i)$); it follows that $\ord \varphi_i = 0 \; \forall i$.

Now consider distinct elements $i(0), \ldots , i(n+1)$ of $\{1, \ldots , m \}$ and minimal projections $e_0 \in M_{r_{i(0)}}, \ldots , e_{n+1} \in M_{r_{i(n+1)}}$, then
\[
\varphi (\be_{i(j)}) = g_{i(j)} \cdot w^*_{i(j)} \be_{i(j)} w_{i (j)} = g_{i (j)} \cdot q_{i (j)}. 
\]
But $g_i$ is nonzero only on $V_i$, and these had strict order less than or equal to $n$, which means that at least two of the sets $V_{i (0)} , \ldots , V_{i (n+1)}$, say $V_{i (j_1)}$ and $V_{i (j_2)}$, do not intersect. Thus $g_{i (j_1)} g_{i (j_2)} = 0$ and $\varphi (\be_{i(j_1)}) \varphi (\be_{i(j_2)}) = 0$, which implies that $\varphi (e_{i(j_1)}) \varphi (e_{i(j_2)}) =0$, hence $\ord \varphi \le n$.
\end{proof}
\en

\bn{\label{bounded-dimension}}
\begin{defn}
A continuous trace $C^*$-algebra $A$ is said to be of locally bounded dimension, if, for every $t \in \hat{A}$, there exists a neighborhood $M$ of $t$ and $r \in \N$ such that each representation $s \in M$ has rank no greater than $r$, i.e.\ $\dim \Hh_s \le r$.
\end{defn}
\en

\bn
\begin{nremarks} \rm
  (i) All homogeneous $C^*$-algebras (which automatically have continuous trace) are, of course, of locally bounded dimension, and so are all subhomogeneous continuous trace $C^*$-algebras.\\
(ii) It is easy to give examples of continuous trace $C^*$-algebras which only have finite-dimensional irreducible representations, but do not have locally bounded dimension.\\
(iii) Since $\hat{A}$ is locally compact for a continuous trace $C^*$-algebra, in the previous definition one may assume $M$ to be compact.
\end{nremarks}
\en

\bn{\label{bounded-rank}}
\begin{lemma}
Let $n,r \in \N$ be given. Then there is $N (n,r) \in \N$ such that the following holds:\\
If $r', k \in \N ,\, r' \le r$, and $\varphi : \C^k \to M_{r'}$ is a c.p.\ map with $\ord \varphi \le n$, then $\varphi (e_i)$ is nonzero for at most $N (n,r)$ of the canonical generators $e_i$ of $\C^k$.
\end{lemma}

\begin{proof}
  After embedding $M_{r'}$ in $M_r$ we may clearly assume $r' = r$ in the assertion. \\
Consider the unit sphere $S^{2r-1} \subset \C^r , \, M_r$ operating on $\C^r$ in the natural way. For every $\xi \in S^{2r-1}$ define
\[
U_{\xi} := \{ \eta \in S^{2r-1} \tei  |\langle \xi | \eta \rangle| > \frac{1}{\sqrt{2}} \} \; ,
\]
then $(U_{\xi})_{\xi \in S^{2r-1}}$ is an open covering of $S^{2r-1}$. Since the sphere is compact, there is a finite subcovering $(V_{\lambda})_{\Lambda}$. Note that if $\eta_1 , \eta_2 \in V_{\lambda}$ ($= U_{\xi}$ for some $\xi$), then
\[
\langle \eta_1 | \eta_2 \rangle = \langle \eta_1 | p_{\xi} | \eta_2 \rangle + \langle \eta_1 | (\be -p_{\xi}) | \eta_2 \rangle \neq 0 \; ,
\]
since
\[
|\langle \eta_1 | p_{\xi} | \eta_2 \rangle| = |\langle \eta_1 | \xi \rangle \langle \xi |\eta_2  \rangle| > \halb
\]
and
\begin{eqnarray*}
  |\langle \eta_1 | (\be-p_{\xi}) | \eta_2 \rangle| & \le & \| (\be - p_{\xi}) \eta_1 \| \| (\be - p_{\xi}) \eta_2\| \\
  & = & \sqrt{\| \eta_1 \|^2 - | \langle \eta_1 | \xi \rangle |^2} \sqrt{\| \eta_2 \|^2 - | \langle \eta_2 | \xi \rangle |^2} \\
& < & 1 - \halb = \halb \; .
\end{eqnarray*}
Here, $p_{\xi}$ denotes the orthogonal projection onto $\C \cdot \xi$ and we used $\langle \mbox{bra} |$-$| \mbox{ket} \rangle$ notation.

Set $N := (n+1) |\Lambda|$ and suppose there are $e_{i_1} , \ldots , e_{i_{N+1}}$ with $\varphi (e_{i_j}) \neq 0 \; \forall\, j$. W.l.o.g.\ we may assume $i_j = j \; \forall \, j$.

Then for each $j$ there is a normed eigenvector $\eta_j$ of $\varphi (e_j)$ w.r.t.\ an eigenvalue $\mu_j > 0$. Now there must be some $\lambda \in \Lambda$ such that at least $n+2$ vectors $\eta_j$ lie in $V_{\lambda}$.

But if $\eta_{j'}, \; \eta_{j''} \in V_{\lambda}$, then 
\[
\langle \eta_{j'} | \varphi (e_{j'}) \varphi (e_{j''}) | \eta_{j''}\rangle = \mu_{j'} \mu_{j''} \langle \eta_{j'} | \eta_{j''} \rangle \neq 0, 
\]
thus $\varphi (e_{j'})$ and $\varphi (e_{j''})$ are not orthogonal.

We therefore obtain a contradiction to $\ord \varphi \le n$.
\end{proof}
\en

\bn{\label{continuous-trace=}}
\begin{theorem}
Let $A$ be a separable continuous trace $C^*$-algebra of locally bounded dimension. Then $\dim \hat{A} \le \cpr A$.
\end{theorem}

\begin{proof}
Suppose for a moment that $A = M_r (\Ch(X))$ for some compact space $X$ and $r \in \N$, so $\Ch(X) \cong \Ch(X) \otimes e_{11} \subset_\her A$. The aim is then to show that $\cpr (\Ch(X)) \le \cpr A$, since by \cite{Wi1}, Proposition 3.18, $\cpr (\Ch(X)) = \dim X$.

Given $\varepsilon > 0$ and $a_1, \ldots , a_l \in \Ch(X) \otimes e_{11}$, one can find a c.p.\ approximation $(F, \psi, \varphi)$ (of $A$) for $a_1, \ldots , a_l$ within $\varepsilon$ and with $\ord \varphi \le \cpr A$. After some extra work one can even assume that $\varphi(\be_F) - e_{11}$ is small, so $(F, \psi|_{\Ch(X) \otimes e_{11}}, \varphi)$ is ``almost'' a c.p. approximation for $\Ch(X)$. However, $\varphi$ maps (the unit ball of) $F$ to $\Ch(X) \otimes e_{11}$ only up to $\varepsilon$; the contributions outside $\Ch(X) \otimes e_{11}$ are certainly small in norm but they might even generate all of $A$ as a $C^*$-algebra. But now Lemma \ref{bounded-rank} (which actually is the key step in the proof) ensures us that this effect cannot be too annoying and that there is $\varphi'' : F \to \Ch(X) \otimes e_{11}$ such that $(F, \psi|_{\Ch(X) \otimes e_{11}}, \varphi'')$ still is a good approximation and such that $\ord \varphi'' \le \ord \varphi \le \cpr A$. 

If $A$ is a continuous trace $C^*$-algebra, the situation is more complicated; the countable sum theorem for covering dimension allows us to ``localize'' the problem, then the idea is again to use Lemma \ref{bounded-rank} to make the above strategy work (the reason why we have to assume locally bounded dimension is that we do not have a version of \ref{bounded-rank} for $\Kh$ instead of $M_r$). 

For each $t \in \hat{A}$ there is a compact neighborhood $M_t$ of $t$, $r_t \in \N$ and $p_t \in A$ with $\|p_t\| \le 1$ such that, for each $s \in M_t$, $p_t (s)$ is a rank-one projection and $\dim \Hh_s \le r_t$. Since $\hat{A}$ is locally compact, there are a compact and an open neighborhood $K_t$ and $V_t$ such that $t \in K_t \subset V_t \subset M_t$. Of course, $(K_t)_{t \in \hat{A}}$ is a covering of $\hat{A}$. But $\hat{A}$ is second countable because $A$ is separable (\cite{Dx}, Proposition 3.3.4), so $(K_t)_{t \in \hat{A}}$ has a countable subcovering $(K_{t_l})_{l \in \N}$. Now if $\dim K_{t_l} \le n := \cpr A$, then by the countable sum theorem for covering dimension we get $\dim \hat{A} \le n$. We thus have to show that $\dim K_t \le n$ for every $K_t$. For convenience, from now on we omit the index $t$. 

We have $\Ch(K) \cong \pi_K(p) (A/J_K) \pi_K(p) \subset_\her A/J_K$; we will simply write $\be_K$ for $\pi_K(p)$. For each $s \in K$ we have $\dim \Hh_s \le r$. Furthermore, $\cpr A/J_K \le n$, since $A/J_K$ is a quotient of $A$. 

Now let $a_1, \ldots , a_l \in \Ch(K)_+$ and $\varepsilon > 0$ be given. Choose some $\delta, \, \eta > 0$ such that
\begin{itemize}
\item[(i)] $2 \cdot N(n,r)^\halb \eta^\frac{1}{16} + 4 \cdot \eta^\frac{1}{8} + 6 \cdot \eta^\frac{1}{4} + \eta < \frac{\varepsilon}{2}$
\item[(ii)] $2 \cdot \delta^\halb N(n,r) + 2 \cdot (3 \cdot \eta^\frac{1}{4} + N(n,r)) \eta^\frac{1}{8} + 2 \cdot (9 + 3 N(n,r))^\halb \eta^\frac{1}{16} < \frac{\varepsilon}{2}$
\item[(iii)] $\frac{\delta^2}{8} > 2 \cdot (3 \cdot \eta^\frac{1}{4} + N(n,r) \eta^\frac{1}{8}) + 2 \cdot (9 + 3 N(n,r))^\halb \eta^\frac{1}{16})$,
\end{itemize}
where $N(n,r)$ comes from Lemma \ref{bounded-rank}.

Let $(F,\psi, \varphi)$ be a c.p.\ approximation (in $A/J_K$) for $\be_K, \, a_1, \ldots, a_l$ within $\eta$ such that $\ord \varphi \le n$. We will modify this approximation in various steps to obtain a c.p.\ approximation of $\Ch(K)$ with the right properties. Set $q:= g_{\eta^\halb} (\psi(\be_K)) \in F$, then $q$ is a projection and 
\[
\varphi(q) \le \frac{1}{\eta^\halb} \varphi \psi (\be_K) \, ;
\]
furthermore
\[
\varphi(q) \ge \varphi \psi (\be_K) - \eta^\halb \cdot \be \, ,
\]
where $\be$ denotes a unit adjoined to $A/J_K$. We then have 
\begin{eqnarray*}
\| \be_K - \varphi(q) \be_K \|^2 & = & \| \be_K (\be - \varphi(q))^2 \be_K \| \\
& \le & \| \be_K (\be - \varphi(q)) \be_K \| \\
& \le & \| \be_K (\be - \varphi \psi(\be_K) + \eta^\halb \cdot \be) \be_K \| \\
& \le & \| \be_K ((1 + \eta^\halb + \eta) \cdot \be - \be_K ) \be_K \| \\
& \le & \eta^\halb + \eta \, ,
\end{eqnarray*}
so 
\[\| \be_K - \varphi(q) \be_K \| < 2 \cdot \eta^\frac{1}{4} \, .
\]
Similarly, we obtain
\begin{eqnarray*}
\| \varphi(q) - \be_K \varphi(q) \|^2 & = & \| (\be - \be_K) \varphi(q)^2 (\be - \be_K) \| \\
& \le &  \| (\be - \be_K) \varphi(q) (\be - \be_K) \| \\
& \le &  \frac{1}{\eta^\halb} \| (\be - \be_K) \varphi \psi(\be_K) (\be - \be_K) \| \\
& < &  \frac{1}{\eta^\halb} \| (\be - \be_K) (\be_K + \eta \cdot \be) (\be - \be_K) \| \\
& \le & \frac{\eta}{\eta^\halb} = \eta^\halb \, ,
\end{eqnarray*}
therefore $\| \varphi(q) - \be_K \varphi(q) \| < \eta^\frac{1}{4}$. As a consequence,
\[
\| \varphi(q) - \be_K \| < 3 \cdot \eta^\frac{1}{4}
\]
and 
\begin{equation}{\label{5}}
\| \varphi(q) - \varphi(q)^2 \| < 4 \cdot \eta^\frac{1}{4} \, .
\end{equation}

We thus obtain for $k = 1, \ldots, l$
\begin{eqnarray*}
\| \varphi(q \psi(a_k) q) - a_k \| & \le & \|  \varphi( q \psi(a_k) q) - \varphi(q) \varphi \psi(a_k) \varphi(q) \| \\
& & + \| \varphi(q) \varphi \psi(a_k) \varphi(q) - \be_K \varphi \psi(a_k) \be_K \| \\
& & + \| \be_K \varphi \psi(a_k) \be_K - \be_K a_k \be_K \| \\
& < & 2 \cdot 2 \cdot \eta^\frac{1}{8} + 6 \cdot \eta^\frac{1}{4} + \eta
\end{eqnarray*}
using \cite{Wi1}, Proposition 1.1.6, and (\ref{5}); for $\be_K$ instead of $a_k$ we have the same estimate.

Misusing our notation, we may write $F$ for $q F q$ (so $q = \be_F$), $\psi$ for $\psi_q$ and $\varphi$ for $\varphi|_{qFq}$ and assume $F \cong M_{r_1} \oplus \ldots \oplus M_{r_s}$. We show that, if $\ord \varphi_i =0$ for some $i$, then either $r_i = 1$ or $\| \varphi_i(\be_{r_i}) \| < \eta^{\frac{1}{8}}$:\\
Suppose $\ord \varphi_i = 0$, $r_i > 1$ and $\| \varphi_i(\be_{r_i}) \| \ge \eta^{\frac{1}{8}}$. Then there is $t \in K$ such that $\| \varphi_i(\be_{r_i})(t) \| \ge \eta^{\frac{1}{8}}$. But since $\ord(\pi_t \verk \varphi_i) =0$ ($\pi_t$ is a representation), each eigenvalue of $\varphi_i(\be_{r_i})(t) \in \Bh(\Hh_t)$ has multiplicity at least $r_i$. It follows that there must be a rank-2 projection $f \in \Bh(\Hh_t)$ such that 
\[
\eta^{\frac{1}{8}} \cdot f \le \varphi_i(\be_{r_i})(t) \le \varphi(\be_F)(t) \, .
\]
But (using $\| \varphi(\be_F)(t) (\be(t) - \be_K(t))\| < \eta^\frac{1}{4}$) one checks that
\[
\| \eta^{\frac{1}{8}} \cdot f - \eta^{\frac{1}{8}} \cdot f \be_K(t) \| < \eta^\frac{1}{8} \, ,
\]
so
\[
\| f - f \be_K(t) \| < \frac{\eta^\frac{1}{8}}{\eta^{\frac{1}{8}}} = 1 \, ,
\]
hence
\[
\| f - f \be_K(t) f \| < 1 \, .
\]
But as $\be_K(t)$ is rank-one, so is $f \be_K(t) f \le f$, a contradiction.

Next set 
\[
p:= \sum_{\{i \tei r_i>1 , \, \ord \varphi_i =0 \}} \be_{r_i} \, ,
\]
then for any $t \in K$ we obtain
\[
\| \varphi(p)(t) \| < N(n,r) \cdot \eta^{\frac{1}{8}} \, :
\]
By Corollary \ref{central-cpr} we have $\ord \varphi \verk \iota \le n$, where $\iota : \C^s \to F$ is the canonical unital embedding. We may thus apply Lemma \ref{bounded-rank} to see that $\varphi(\be_{r_i})(t)$ is nonzero for at most $N(n,r)$ values of $i$, and by the preceding observation, $\| \varphi_i(\be_{r_i}) \| < \eta^\frac{1}{8}$ for those $i$.

We now have for $k =1, \ldots ,l$
\begin{eqnarray*}
\| \varphi((\be_F - p) \psi(a_k) (\be_F - p)) - a_k \| & < & 2 \cdot N(n,r)^\halb \cdot \eta^{\frac{1}{16}} \\
& & + 2 \cdot 2 \cdot \eta^\frac{1}{8} + 6 \cdot \eta^\frac{1}{4} + \eta \\ 
& < & \frac{\varepsilon}{2} \, ; 
\end{eqnarray*}
the same estimate holds for $\be_K$ instead of $a_k$.

So, oncemore misusing our notation, we write $F$ for $(\be_F - p) F (\be_F - p)$, $\psi$ for $\psi_{\be - p}$ and $\varphi$ for $\varphi|_{(\be_F - p) F (\be_F - p)}$. Also, we still write $F = M_{r_1} \oplus \ldots \oplus M_{r_s}$; this should not cause confusion.

After these modifications, we now have a c.p.\ approximation $(F, \psi, \varphi)$ (of $A/J_K$) for $\be_K, \, a_1, \ldots, a_k$ within $\frac{\varepsilon}{2}$ with the following properties:
\begin{itemize}
\item[(a)] $\| \varphi(\be_F) - \be_K \| < 3 \cdot \eta^\frac{1}{4} + N(n,r) \cdot \eta^\frac{1}{8}$
\item[(b)] $\ord \varphi \le n$
\item[(c)] $\ord \varphi_i = 0 \mbox{ iff } M_{r_i} = \C$.
\end{itemize}

For $x \in F_+$ with $\|x\| \le 1$, (a) in particular implies
\def\theequation{$\ast \ast$}
\begin{equation}{\label{6}}
\| \varphi(\be_F) \varphi(x) - \varphi(x)\| < (9 + 3 \cdot N(n,r))^\halb \cdot \eta^\frac{1}{16}
\end{equation}
\def\theequation{$\ast$}
because (again by \cite{Wi1}, Proposition 1.1.6) 
\begin{eqnarray*}
\| \varphi(\be_F) - \varphi(\be_F)^2 \| & \le & \| \varphi(\be_F) - \be_K \| + \| \be_K - \varphi(\be_F)^2 \| \\
& < & 3 \cdot ((3 + N(n,r))) \cdot \eta^\frac{1}{8} \, .
\end{eqnarray*}
 
Now define $\varphi':F \to \Ch(K) \subset_\her A/J_K$ by
\[
\varphi'(\,.\,) := \be_K \varphi(\,.\,) \be_K \, ,
\]
then by a) and (\ref{6}) we have for $x \in F_+$ with $\|x\| \le 1$ 
\begin{eqnarray*}
\| \varphi'(x) - \varphi(x)\| & \le & \| \be_K \varphi(x) \be_K - \varphi(\be_F) \varphi(x) \varphi(\be_F) \| \\
& & + \| \varphi(\be_F) \varphi(x) \varphi(\be_F) - \varphi(x) \| \\
& < &  2 \cdot (3 \cdot \eta^\frac{1}{4} + N(n,r) \cdot \eta^\frac{1}{8}) + 2 \cdot (9 + 3 N(n,r))^\halb \eta^\frac{1}{16} \,.
\end{eqnarray*}
Next define $\varphi'' : F \to \Ch(K)$ by setting
\[
\varphi_i''(\,.\,) := g_{\frac{\delta}{2}, \delta} (\varphi'(\be_{r_i})) \varphi_i'(\,.\,) g_{\frac{\delta}{2}, \delta} (\varphi'(\be_{r_i})) \, .
\]
For $0 \le x_i \in M_{r_i}$, $\|x_i\| \le 1$,
\begin{eqnarray*}
\| (\be - g_i) \varphi_i'(x_i) \|^2 & = & \| (\be - g_i) \varphi_i'(x_i)^2 (\be - g_i) \| \\
& \le & \| (\be - g_i) \varphi_i'(x_i) (\be - g_i) \| \\
& \le & \| (\be - g_i) \varphi_i'(\be_{r_i}) (\be - g_i) \| \\
& < & \delta \, ,
\end{eqnarray*}
where we have written $g_i$ for $g_{\frac{\delta}{2}, \delta} (\varphi'(\be_{r_i}))$. As a consequence,
\[
\| \varphi_i''(x_i) - \varphi_i'(x_i) \| < 2 \cdot \delta^\halb \, .
\]
Furthermore, again from Lemma \ref{bounded-rank} and Corollary \ref{central-cpr} we know that for each $t \in K$
\[
\varphi_i'(\be_{r_i})(t) = \be_K(t) \varphi_i(\be_{r_i})(t) \be_K(t)
\]
is nonzero for at most $N(n,r)$ values of $i$; of course $\varphi_i''(\be_{r_i})(t)$ is zero if $\varphi_i'(\be_{r_i})(t)$ is. So we have that for all $x \in F_+$ with $\|x\| \le 1$
\begin{eqnarray*}
\| \varphi''(x) - \varphi'(x) \| & = & \sup_t \| \varphi''(x)(t) - \varphi'(x)(t) \| \\
& \le & \sup_t ( \sum_i \| \varphi_i''(x_i) (t) - \varphi_i'(x_i) (t) \| ) \\
& < & 2 \cdot \delta^\halb N(n,r) \, .
\end{eqnarray*}
Thus
\begin{eqnarray*}
\| \varphi''(x) - \varphi(x) \| & < & 2 \cdot \delta^\halb N(n,r) + 2 \cdot (3 \eta^\frac{1}{4} + N(n,r) \eta^\frac{1}{8}) \\
& & + 2 \cdot (9 + 3 N(n,r))^\halb \eta^\frac{1}{16}) \\
& \le & \frac{\varepsilon}{2} \\
\end{eqnarray*}
and
\begin{eqnarray*}
\| \varphi'' \psi (a_k) - a_k \| & \le & \| \varphi'' \psi (a_k) - \varphi \psi (a_k) \| \\
& & + \| \varphi \psi(a_k) - a_k \| \\
& < & \frac{\varepsilon}{2} + \frac{\varepsilon}{2} = \varepsilon \, .
\end{eqnarray*}

$\ord \varphi'' \le n$:

If $\varphi_i''(\be_{r_i}) \varphi_{\bi}''(\be_{r_{\bi}}) \neq 0$ for some $i \neq \bi \in \{1, \ldots, s\}$, then there is $t \in K$ such that $0 \neq \varphi_i''(\be_{r_i})(t)$, $\varphi_{\bi}''(\be_{r_\bi})(t) \in \C$. But then $\varphi_i'(\be_{r_i})(t)$, $\varphi_\bi'(\be_{r_\bi})(t) \ge \frac{\delta}{2}$, i.e. $\| \varphi_i'(\be_{r_i}) \varphi_\bi'(\be_{r_\bi}) \| \ge \frac{\delta^2}{4}$ (here we have used that $\varphi', \, \varphi'' : F \to \Ch(K) \subset A/J_K$ and that each $t \in K$ represents an irreducible, hence one-dimensional representation of $\Ch(K)$). As a consequence,
\begin{eqnarray*}
\| \varphi_i(\be_{r_i}) \varphi_\bi(\be_{r_\bi}) \| & \ge & \left| \| \varphi_i(\be_{r_i}) \varphi_\bi(\be_{r_\bi}) - \varphi_i'(\be_{r_i}) \varphi_\bi'(\be_{r_\bi}) \| - \| \varphi_i'(\be_{r_i}) \varphi_\bi'(\be_{r_\bi}) \| \right| \\
& \ge & \frac{\delta^2}{4} - 2 \cdot ( 2 \cdot (3 \cdot \eta^\frac{1}{4} + N(n,r) \eta^\frac{1}{8})\\
& & + 2 \cdot (9 + 2 N(n,r))^\halb \eta^\frac{1}{16} ) \\
& > & 0 \, .
\end{eqnarray*}
Of course $\ord \varphi_i =0$ implies $\ord \varphi_i'' =0$, since in this case $r_i =1$. We can thus apply Proposition \ref{dominated-order} and obtain $\ord \varphi'' \le \ord \varphi \le n$.

So we have constructed a c.p.\ approximation $(F, \psi|_{\Ch(K)}, \varphi'')$ (of $\Ch(K)$) for $a_1, \ldots, a_l$ within $\varepsilon$ and with $\ord \varphi'' \le n$, so $\cpr \Ch(K) \le n$. Now Proposition 3.18 of \cite{Wi1} says that $\dim K \le n$ and our proof is complete.
\end{proof}
\en

\bn
\begin{nexamples}
(i) For any separable locally compact space $X$ and $r \in \N , \; \Ch_0 (X) \otimes M_r$ is $r$-homogeneous, so $\cpr (\Ch_0 (X) \otimes M_r) = \dim X$.\\
(ii) The rational rotation algebras $A_{\theta}$ are homogeneous with spectrum $\T^2$, so $\cpr (A_{\theta}) = 2$ for $\theta$ rational.
\end{nexamples}
\en

%\input{sec7}
%\input{sec8}
%\newpage
%\addtocontents{toc}{{\sc References}}

%\newpage
\addcontentsline{toc}{section}{\sc Contents}
\tableofcontents
%\input{test}

%\clearpage
%\thispagestyle{empty}
%\vspace*{\fill}

%\clearpage
%\thispagestyle{empty}
%\vspace*{\fill}

%\clearpage
%\thispagestyle{empty}
%\vspace*{\fill}

\end{document}